\documentclass[10pt]{article}

\usepackage{enumerate}
\usepackage{graphicx}
\usepackage{amsmath}
\usepackage{amssymb}

\usepackage{amsthm}%entornos theorem
\usepackage{amscd}
\usepackage[latin1]{inputenc}
\def\r{\mathbb{R}}
\def\b{\mathbb{B}}
\def\d{\mathbb{D}}
\def\n{\mathbb{N}}
\def\c{\mathbb{C}}
\def\s{\mathbb{S}}

\def\nor{{\cal N}}
\def\pro{{\cal P}}
\newcommand{\df}{ \stackrel{\rm def}{=}}
\newcommand{\cte}{\text{\rm const}}
\newcommand{\rth}{\r^3}
\newcommand{\metri}[1]{{\cal S}_{#1}}

\newtheorem{teorema}{Theorem}

\newtheorem*{teoremaintroA}{Theorem A}
\newtheorem*{teoremaintroB}{Theorem B}
\newtheorem*{teoremaintroC}{Theorem C}

\newtheorem{claim}{Claim}[section]
\newtheorem{lema}{Lemma}

\newtheorem{proposicion}{Proposition}
\newtheorem{remark}{Remark}

\newcommand{\dist}{\operatorname{dist}}
\newcommand{\re}{\operatorname{Re}}

\newcommand{\met}[1]{\operatorname{\cal S}_{#1}}

\newcommand{\intc}{\operatorname{Int}}

\newcommand{\longui}{\operatorname{\ell}}

\newfont{\goti}{cmr10 at 11pt}
\begin{document}

\title{ Complete proper minimal surfaces in convex bodies of $\r^3$ (II): The behavior of the limit set}
\author{\\ Francisco Martín\thanks{Research partially supported by MEC-FEDER Grant no. MTM2004 - 00160.} \and \\ Santiago Morales$^*$}
\date{May 23, 2005}
\maketitle
\begin{abstract} 

Let $D$ be a regular strictly convex  bounded domain of $\r^3$, and consider a regular Jordan curve $\Gamma \subset \partial D$. Then, for each $\varepsilon>0$, we obtain the existence of a complete proper minimal immersion $\psi_\varepsilon :\d \rightarrow D$ satisfying that the Hausdorff distance $\delta^H(\psi_\varepsilon(\partial \d), \Gamma) < \varepsilon,$ where $\psi_\varepsilon(\partial \d)$ represents the limit set of the minimal disk $\psi_\varepsilon(\d).$

This result has some interesting consequences. Among other things, we can prove that any  bounded regular domain $R$ in $\r^3$  admits a complete proper minimal immersion $\psi: \d \longrightarrow R$.
\vskip .5cm

\noindent {\em 2000 Mathematics Subject Classification. Primary 53A10; Secondary 49Q05, 49Q10, 53C42.
}
\end{abstract}
\section{Introduction and background} \label{sec:intro}

Last few years have seen an important progress on many long-standing problems in global theory of complete minimal surfaces in $\mathbb{R}^3$. One of these  has been the Calabi-Yau problem, which dates back to the 1960s.  Calabi asked wheter or not it is possible for a complete minimal surface in $\r^3$ to be contained in the ball $\b=\{x \in \r^3 \, | \, \|x\|<1\}.$ Much work has been done on it over the past four decades. The most important result in this line was obtained by N. Nadirashvili in \cite{nadi} where he constructed a complete minimal surface in $\b$. After Nadirashvili's negative solution to Calabi's question, the conjecture was revisited by S.-T. Yau in \cite{yau}, where he stated new questions related to the embeddedness and properness of surfaces of this type.

Regarding the existence of complete embedded minimal surfaces in a ball, T. Colding and W. Minicozzi \cite{c-m} have proved that a complete embedded minimal surface with finite topology in $\r^3$ must be properly embedded in $\r^3$. In particular it cannot be contained in a ball. Very recently, Colding-Minicozzi result has been generalized  in two different directions. On one hand W. H. Meeks III, J. Pérez and A. Ros \cite{mpr-2} have proved  that if $M$ is a complete embedded minimal surface in $\r^3$ with finite genus and a countable number of ends, then $M$ is properly embedded in $\rth$. On the other hand, Meeks and Rosenberg \cite{mr-lam} have obtained that if a complete embedded minimal surface $M$ has injectivity radius $I_M >0$,  then $M$ is  proper in space.

Concerning properness, it is important to note that Nadirashvili's technique did not guarantee the immersion was 
proper. In  \cite{tran}, F. Martín and S. Morales  introduced an additional ingredient into Nadirashvili's machinery in order to
produce a complete minimal disk which is properly immersed in a ball of $\r^3$.  Recently \cite{convex-I}, they  improved on their original techniques and were able to show that every  convex domain (not necessarily bounded or smooth) admits a complete properly immersed minimal disk.

The present paper can be considered as a continuation of the above mentioned work developed by the authors about the construction of complete proper minimal surfaces in (open) convex bodies of Euclidean space. After the discovering of those examples a natural question arose: What is the asymptotic behavior of such a  surface ? If we consider a proper minimal immersion $\psi: \d \rightarrow C$, where $C$ is an open convex body, then we define the limit set as $\psi(\partial \d) \df \overline{\psi(\d)} -\psi(\d)$. It is obvious that $\psi$ is proper if, and only if, $\psi(\partial \d) \subset \partial C.$ Furthermore, one can easily check that $\psi(\d)$ is closed and connected. In this paper we show:
\begin{teoremaintroA}
Let $C$ be a regular strictly convex\footnote{{\bf strictly convex} means that the principal curvatures of $\partial C$ associated to the inward pointing unit normal are positive everywhere.}  bounded domain of $\r^3$, and consider a regular Jordan curve $\Gamma \subset \partial C$. Then, for each $\varepsilon>0$, we obtain the existence of a complete proper minimal immersion $\psi_{(\Gamma,\varepsilon)} :\d \rightarrow C$ satisfying that the Hausdorff distance $\delta^H(\psi_{(\Gamma,\varepsilon)}(\partial \d), \Gamma) < \varepsilon,$ where $\psi_{(\Gamma,\varepsilon)}(\partial \d)$ represents the limit set of the minimal disk $\psi_{(\Gamma,\varepsilon)}(\d).$
\end{teoremaintroA}
The main obstacle in the study of the asymptotic behavior of a complete proper minimal surface in a convex domain is that all the methods of construction up until today are implicit. So, it was almost impossible to obtain any control about the behavior of the known examples near their ends. From this point of view, the new methods of construction introduced in this paper are significant. The main of these tools is Theorem \ref{th:meeks} whose proof is based on a Meeks' idea that  appeared first in \cite{convex-I} but that has been entirely developed and exploited in this article.

\begin{figure}[htbp]
	\begin{center}
		\includegraphics[width=0.5\textwidth]{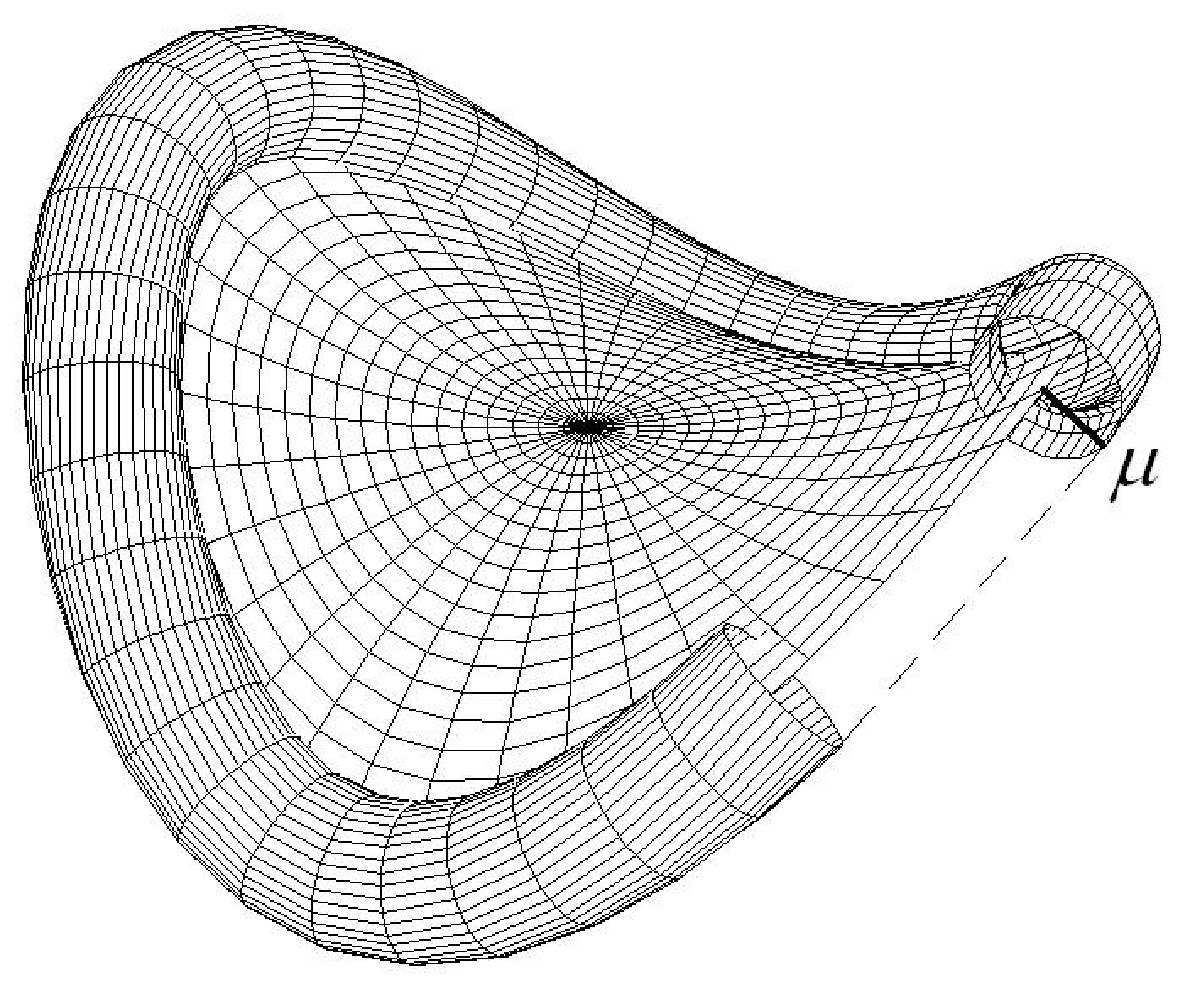}
	\end{center}
	\caption{The minimal disk $X(P)$ and the tube $T(\Gamma,\mu).$}
	\label{fig:disco-tubo}
\end{figure}

Theorem \ref{th:meeks} represents by itself an interesting density result. It asserts that any minimal disk with boundary $D$ can be approximated (in terms of uniform convergence) by a complete minimal disk $\widetilde D.$ Moreover, it is possible to find a thin tube around $\partial D$ such that the part of $\widetilde D$ which lies in the exterior of this tube is compact. In other words, the part of $\widetilde D$ where the intrinsic metric exploits is contained in the interior of this thin tube (see Figure \ref{fig:disco-tubo}.)

The other important ingredient in the proof of the main theorem is Lemma \ref{lem:propia}. This approximation lemma is crucial to obtain the properness in Theorem A and it was proved in \cite{convex-I}. It essentially asserts that a minimal disk with boundary can be perturbed  outside a compact set in such a way that the boundary of the resulting surface achieves the boundary of a prescribed convex domain.  Hypothesis of $C$ being strictly convex is crucial at this point, otherwise we could not obtain an upper bound  for the distances between the boundary of a minimal disk and the boundary of the deformed one.

It is natural to ask what is the limit of the complete proper minimal surfaces given by Theorem A as $\varepsilon$ tends to $0$. We would like to point out that the limit as $\varepsilon$ goes to zero of $\psi_{(\Gamma, \varepsilon)}$ exists, but it is not complete. Actually, this limit coincides with the minimal disk (with boundary) spanned by the curve $\Gamma.$

Theorem A has some interesting consequences. We would like to point out two of them.
\begin{teoremaintroB}
Let $C$ be a regular, strictly convex  bounded domain, and consider a  {\em connected compact set} $K \subset \partial C$. Then, for each $\varepsilon>0$, there exists a complete proper minimal immersion $\varphi_{(K,\varepsilon)}:\d \rightarrow C$ satisfying that the Hausdorff distance $\delta^H(\psi_{(K,\varepsilon)}(\partial \d), K) < \varepsilon.$
\end{teoremaintroB}
The above theorem  follows from the fact that Jordan curves are dense in the space of compact sets of $\partial C$ with the Hausdorff metric. Among other things, Theorem B says to us that the limit set of a complete proper minimal surface can be very small. This means that we can work with small pieces of the boundary of a given domain in orther to prove the following: 
\begin{teoremaintroC}
Every bounded  domain with regular boundary admits a complete properly immersed minimal disk.
\end{teoremaintroC}
In contrast to these existence results for complete properly immersed minimal disks in bounded domains, Meeks, Nadirashvili and the first author \cite{mmn} have constructed domains of  $\rth$ which do not contain any complete proper minimal surface with finite topology. It is our belief that these   open domains are in fact {\em universal} according to the following definition: A connected region of space which is open or the closure of an open set is {\it universal for minimal surfaces}, if every complete properly immersed minimal surface in the region is recurrent for Brownian motions. In particular, a bounded domain is universal if and only if it contains no complete properly immersed minimal surfaces.
\begin{figure}[bpht]
	\begin{center}
		\includegraphics[width=0.4\textwidth]{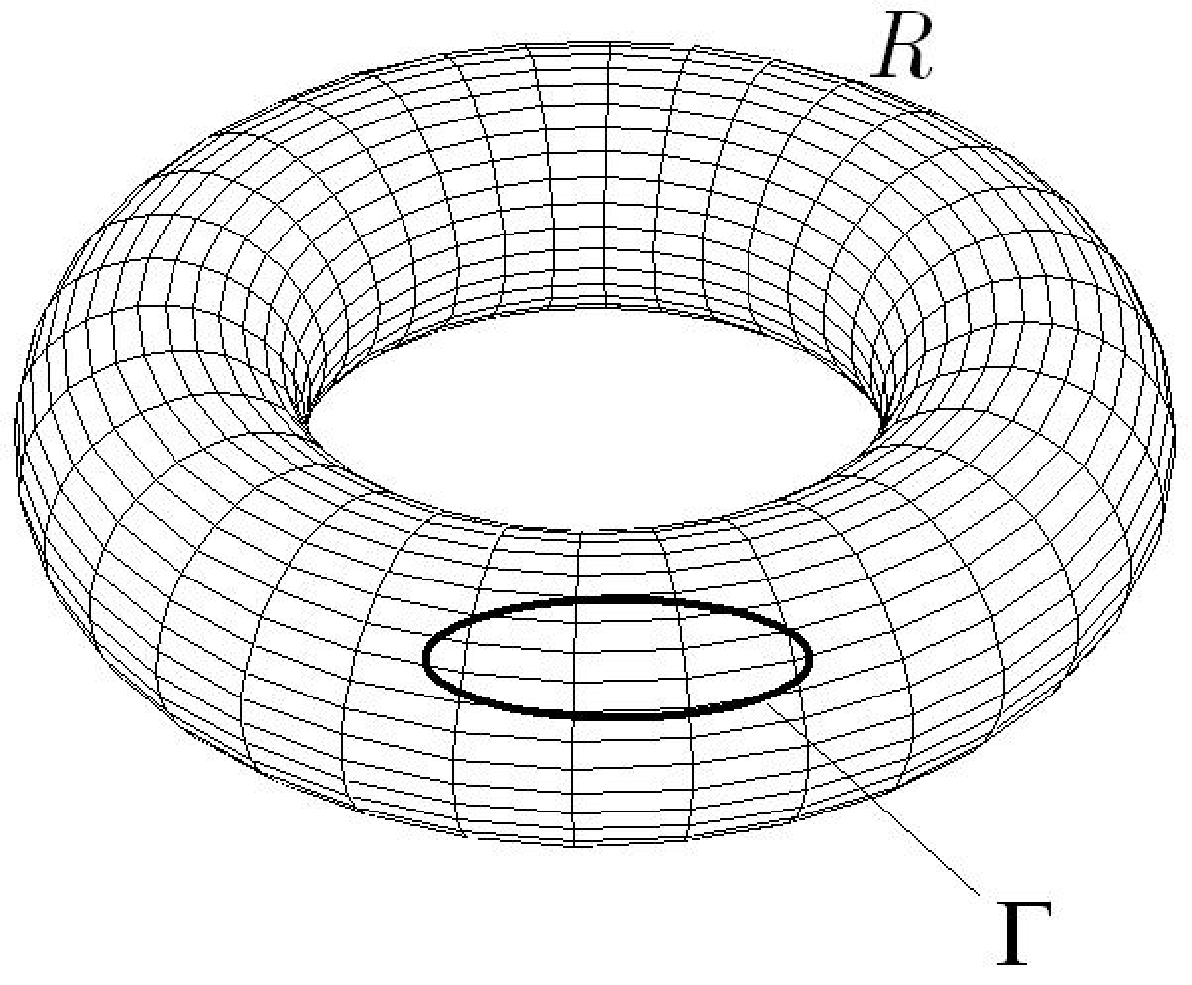}
	\end{center}
	\caption{Theorem C states that every  bounded domain with regular boundary admits  a  complete properly immersed minimal surface. Besides,  the boundary limit set is close to a small simple closed curve $\Gamma$ on the positively curved part of the boundary of the domain.}
	\label{fig:torus}
\end{figure}

As we mentioned before, Colding and Minicozzi proved that any complete embedded minimal surface in $\rth$ with finite topology is properly embedded in $\rth$. By results of Meeks and Rosenberg, \cite{mr8, mr6}, any properly embedded minimal surface of finite topology in $\rth$ is recurrent for Brownian motion. Hence, every domain in $\rth$ is universal for embedded minimal surfaces of finite topology. Finally, we remark that Collin, Kusner, Meeks and Rosenberg \cite{ckmr} proved that any properly immersed minimal surface with boundary in a closed convex domain in $\rth$ has full harmonic measure on its boundary. 

The paper is displayed as follows. Section \ref{sec:lemmata} is dedicated to state and demonstrate the preliminary results that we need to obtain the main theorems. These main theorems are proved in Section \ref{sec:main} where we also establish  some results which show that Theorem A is sharp. Thus, at the end of Section \ref{sec:main} we show that Theorem A is false if the boundary of the domain contains an open set where the mean curvature is non positive.

\subsection{ Minimal surface background}

Given $X=(X_1,X_2,X_3):M \longrightarrow \r^3$ a conformal minimal immersion we denote by $g:M \longrightarrow \overline{\c}=\c \cup \{\infty \}$ its stereographically projected Gauss map that is a meromorphic function and by $\phi_3$ the holomorphic differential defined as $\phi_3=d X_3+ \star \, {\rm i} \, d X_3$, where $\star$ denotes the Hodge operator on $M$. The pair $(g,\phi_3)$ is usually referred to as the Weierstrass data of the minimal surface, and the minimal immersion $X$ can be expressed, up to translations, solely in terms of these data as 
\begin{equation}\label{eq:inmersion} \hspace{-0.5in} X={\rm Re}\int^z (\phi_1,\phi_2,\phi_3)={\rm Re}\int^z \left(\frac{1}{2} \left( \frac{1}{g}-g \right),\frac{{\rm i}}{2}\left(\frac{1}{g}+g \right),1 \right) \phi_3 \;, \end{equation}
where ${\rm Re}$ stands for real part and $z$ is a conformal parameter on $M$. The pair $(g,\phi_3)$ satisfies certain compatibility conditions:
\begin{equation} \label{eq:conforme}
\phi_1^2+\phi_2^2+\phi_3^2=0;
\end{equation}
\begin{equation} \label{eq:bilbao}
\|\phi_1\|^2+\|\phi_2\|^2+\|\phi_3\|^2 \neq 0;
\end{equation}
and all periods of the $\phi_j$ are purely imaginary; $j=1,2,3$.

Conversely, if $M$ is a Riemann surface, $g:M \to \overline{\c}$ is a meromorphic function and $\phi_3$ is a holomorphic one-form on $M$ fulfilling the conditions \eqref{eq:conforme} and \eqref{eq:bilbao} then the map $X:M \to \r^3$ given by
\eqref{eq:inmersion} is a conformal minimal immersion with Weierstrass data $(g,\phi_3)$.

Condition ii) stated above deals with the independence of  \eqref{eq:inmersion} on the integration path, and it is usually  called the period problem.
 In this article, all the minimal immersions are defined on simply connected domains of $\c$. Then, the Weierstrass 1-forms have no periods, and so the only requirements are (\ref{eq:conforme}) and (\ref{eq:bilbao}). In this case, the differential $\eta \df \phi_3/g$ can be written as $\eta=f(z) \, dz$. 
The metric of $X$ can be expressed as
\begin{equation} \label{eq:metric}
{\metri X}^2=\tfrac12 \|\phi\|^2=\left(\tfrac12\left(1+|g|^2\right) \; |f| \; |dz|\right)^2.
\end{equation}
Throughout the paper, we will use several orthonormal bases of $\r^3$. Given $S$ an orthonormal basis and $v \in \r^3$, let $v_{(k,S) }$ denote the $k$-th coordinate of $v$ in $S$. The first two coordinates of $v$ in this basis will be represented by $v_{(*,S)}=\left(v_{(1,S)},v_{(2, S)} \right)$.

Given a curve $\alpha$ in $\Omega$, by $\longui(\alpha,X)$ we mean the length of  $\alpha$ with respect to the metric $\metri X$. Given a subset $W \subset \Omega$, we define:
\begin{itemize}
\item $\dist_{(W,X)}(p,q)=\inf \{\longui(\alpha, \metri X) \: | \: \alpha:[0,1]\rightarrow W, \; \alpha(0)=p,\alpha(1)=q \}$, for any $p,q\in W$;
\item $\dist_{(W,X)}(T_1,T_2)=\inf \{\dist_{(W,\metri X)}(p,q) \;|\;p \in T_1, \;q \in T_2 \}$, for any $T_1, T_2 \subset  W$;
\end{itemize}
The Euclidean metric on $\c$ will be denoted as $\langle \cdot,\cdot \rangle$. Note that ${\metri X}^2=\lambda_X^2 \,\langle \cdot,\cdot \rangle$, where the conformal coefficient $\lambda_X$ is given by (\ref{eq:metric}). 
 
Given a domain $D\subset \c$, we will say that a function, or a 1-form, is harmonic, holomorphic, meromophic,... on $\overline D$, if it is harmonic, holomorphic, meromorphic,... on a domain containing $\overline D$.

Let $P$ be a simple closed polygonal curve in $\c$. We let $\intc P$
denote the bounded connected component of $\c - P.$ We will assume that the origin is in the interior region of all the polygons that appears in the paper.  Given
$\xi>0$, small enough, we define $P^\xi$ to be the parallel polygonal
curve in $\intc P$, satisfying the property that the distance between parallel sides is equal to $\xi$. Whenever we write $P^\xi$ in the paper we are assuming that $\xi$ is small enough to  define the polygon properly. 

\subsection{Background on convex bodies and Hausdorff distance}
 
Given $E$ a bounded regular convex domain of $\r^3$ and $p \in
\partial E$, we will let $\kappa_2(p) \geq \kappa_1(p) \geq 0$ denote the 
principal curvatures of $\partial E$ at $p$ (associated to the inward
pointing unit normal.) Moreover, we write: $$\kappa_1(\partial E) \df \mbox{min} \{ \kappa_1(p) \: : \: p \in \partial E \} , \qquad \kappa_2(\partial E) \df \mbox{max} \{ \kappa_2(p) \: : \: p \in \partial E \}.$$ 
If we consider $\nor: \partial E \rightarrow \s^2$ the outward
pointing unit normal or Gauss map of $\partial E$, then there exists a constant $a>0$ (depending on $E$) such that $\partial E_t=\{p+ t\cdot \nor(p) \; : \; p \in \partial E\}$ is a regular (convex) surface $\forall t \in [-a, +\infty[$. Let  $E_t$ denote the convex domain bounded by $\partial E_t$. 
The normal projection to $E$ is represented as 
$$\pro_E:\r^3\setminus E_{-a}\longrightarrow \partial E,$$
$$p+t \cdot \nor(p) \mapsto p.$$

For a curve $\Upsilon$ in $\r^3$ and a  real $r>0$, we define the tube of radius $r$ along $\Upsilon$ in the following way:
$$T(\Upsilon,r) = \Upsilon + \b(0,r),$$
where $\b(0,r)=\{ p \in \r^3 \; : \; \|p\|<r\}.$

The set ${\cal C}^n$ of convex bodies  of $\r^n$ can be made into a metric space in several geometrically reasonable ways. The Hausdorff metric is particularly convenient and applicable. The natural domain for this metric is the set ${\cal K}^n$ of the nonempty compact subsets of $\r^n$. 

For $C$, $D\in {\cal K}^n$ the {\em Hausdorff distance} is defined by:
$$\delta^H(C,D)=\max \left\{ \sup_{x \in C} \inf_{y \in D} \|x-y \|, \sup_{y \in D} \inf_{x \in C} \|x-y \| \right\}$$
or, equivalently, by
$$ \delta^H(C,D)=\min \left\{ \lambda \geq 0 \; | \; C \subset D+ \lambda \b^n, \; D \subset C+ \lambda \b^n \right\},$$
where $\b^n=\{p \in \r^n \; | \; \|p \| <1\}. $ Then $\delta^H$ is a metric on ${\cal K}^n$, the {\em Hausdorff metric.} For more details we refer to \cite{sch}.
\section{Preliminary Lemmas} \label{sec:lemmata}
As we indicated in the introduction, the proofs of the main theorems of the paper require the technical results of this section. To be more precise, the two principal tools are Theorem \ref{th:meeks} and Lemma \ref{lem:propia}. Lemma \ref{lem:nadi} is a suitable combination of the two previous  results that allows us to prove Theorem \ref{th:mari}. Regarding Lemma  \ref{lem:complete}, it is a necessary instrument in the demonstration of Theorem \ref{th:meeks}.

\begin{lema}[Completeness Lemma] \label{lem:complete}
Consider $\Gamma$ a closed analytic  curve.
Let  $P$ be a polygon, $X:\overline{\intc P} \rightarrow \r^3$ a conformal minimal immersion, and $r,\epsilon$ positive constants satisfying: \begin{enumerate}[\rm 1.] 
\item $X(\overline{\intc P - \intc P^\epsilon } ) \subset T(\Gamma, r);$
\item $X(P^\epsilon)$ contains a cycle which is homologically equivalent to $\Gamma$ in $T(\Gamma, r).$ 
\end{enumerate}

\noindent Then, for all $s>0$ there exist a polygon $\widetilde P$ and a conformal minimal immersion $\widetilde X : \overline{\intc \widetilde P} \rightarrow  \r^3$ verifying:
\begin{enumerate}[\rm ({a.}1)]
\item $\intc P^\epsilon \subset \intc \widetilde P \subset \overline{\intc \widetilde P} \subset \intc P;$
\item $\|X(z)-\widetilde X(z) \| < \epsilon , $ $\forall z \in \intc P^\epsilon;$
\item $\dist_{(\overline{\intc \widetilde P},\widetilde X)}(z, P^\epsilon)> s,$ $\forall z \in \widetilde P;$
\item $\widetilde X\left( \overline{\intc \widetilde P - \intc P^\epsilon} \right) \subset T( \Gamma,R)$ where $R=\sqrt{(2 s)^2+r^2}+\epsilon;$
\item $X(P^\epsilon)$ and $\widetilde X(P^\epsilon)$ are homologically equivalent in $T(\Gamma, R)$;
\item   $\widetilde X (\widetilde P)$ contains a cycle with the same homology than $\Gamma$ in the tube $T(\Gamma, R)$.

\end{enumerate}
\end{lema}
\vskip .3cm

As we mentioned before, Lemma \ref{lem:complete} is merely a tool in the proof of the next theorem, that we have called Meeks' trick because it is based on an idea that W.H. Meeks III suggested to us in 2004. Roughly speaking this theorem asserts that complete minimal disks are ``dense'' in the space of minimal disks with  boundary. 

\begin{teorema}[Meeks' trick] \label{th:meeks}
Let $U \subset \c$ be a bounded domain, and $P \subset U$ a polygon. Consider $X:U \rightarrow \r^3$ a conformal minimal immersion, with $X(0)=0.$
Then, for every $\mu>0$ there exists a simply connected domain $\Sigma \subset \c$ and a complete minimal immersion $\widehat X:\Sigma \rightarrow \r^3$, with $\widehat X(0)=0$ such that:
{\samepage \begin{enumerate}[\rm ({b.}1)]
\item $\overline{\intc P} \subset \Sigma \subset \overline \Sigma \subset U;$
\item $\|X(z)-\widehat X(z) \|< \mu,$ $\forall z \in \intc P;$
\item $\widehat X(\Sigma - \intc P) \subset T(X(P),\mu);$
\item $\widehat X(P)$ is homologous to $X(P)$ in the open neighborhood $T(X(P),\mu).$
\end{enumerate}}
\end{teorema}

The next lemma  was obtained in \cite{convex-I} (see Lemma 1 and Remark 3). We have stated it here just to make this paper self-contained. We would like to point out that the assumption of $C$ being strictly convex is essential in this lemma, otherwise item (c.5) has no sense. 

\begin{lema}[Properness Lemma, \cite{convex-I}] \label{lem:propia}
Let $E$ and $E'$ be two  regular bounded strictly convex domains in $\r^3$, with $0\in E \subset \overline{E} \subset E'$.  Let $X:O\longrightarrow\r^3$ be a conformal minimal immersion defined on a simply connected domain $O$, $0\in O$, with $X(0)=0$. Consider a polygon $ P$ with $ P\subset  O$,  satisfying:
\begin{equation}\label{eslabon}
 X(O\setminus \intc P)\subset E ' - \overline{E}.
\end{equation}
Then, for any $ b_1, b_2>0$, such that $E'_{-b_2}$ and $E_{-2 b_2}$ exist, there exist a polygon $Q$ and a conformal minimal immersion $Y:\overline{\intc Q}\longrightarrow \r^3$, with $Y(0)=0$, such that:
\begin{enumerate}[\rm ({c.}1)]
\item $ P\subset \intc Q\subset\overline{\intc Q}\subset  O$;
\item $\|Y(z)- X(z)\|< b_1$, $\forall z\in\overline{\intc  P}$;
\item $Y(Q)\subset E'- E'_{- b_2}$;
\item $Y(\intc Q\setminus\intc  P)\subset \r^3\setminus E_{-2 \, b_2}$;
\item For any $z\in \intc Q\setminus\intc  P$, one has $\|X(z)-Y(z) \| < M(b_1,b_2,E,E')$, where 
$$M(b_1,b_2,E,E')=\left[\left( 2 b_2+\frac1{\kappa_1(\partial E)}+ \delta^H(E,E')\right)^2 -\left(2 b_2+\frac1{\kappa_1(\partial E)}\right)^2+(2 b_2)^2\right]^{\frac 12}+b_1.$$
\item $X(P)$ and $Y(P)$ are homologous cycles in $T \left(X(P), M(b_1,b_2,E,E') \right).$
\end{enumerate}
\end{lema}
The last lemma of this section is a combination of Theorem \ref{th:meeks} and Lemma \ref{lem:propia}. The proof of our main result (Theorem \ref{th:mari}) consists of constructing a sequence of minimal disks with boundary whose limit is the immersion which proves the theorem. The sequence is defined in a recursive way and the tool to obtain an element of that sequence from the previous one in Lemma \ref{lem:nadi} below.
\begin{lema} \label{lem:nadi}
Let $\Gamma$ be a smooth  Jordan curve of $\r^3$ and let $\rho>0$ be a sufficiently small constant so that the tube $T(\Gamma, \rho)$ is  homeomorphic to a solid torus. Let $E$ and $E'$ be two regular bounded strictly convex domains in $\r^3$, with $0\in E \subset \overline{E} \subset E'$. Consider $P$  a polygon, $X: \overline{\intc P} \longrightarrow\r^3$ be a conformal minimal immersion, with $X(0)=0$, and $\varepsilon$, $a$, $b$ and $c$ positive constants, such that:
\begin{enumerate} 
\item $X(\overline{\intc P\setminus\intc P^\varepsilon} )\subset E\setminus E_{-a};$
\item $\delta^H(\Gamma, X(\overline{\intc P - \intc P^\varepsilon}) \leq c;$
\item $X(P)$ is homologous to $\Gamma$ in $T(\Gamma,c);$
\item $c+M(a,b,\varepsilon, E,E')< \rho,$ where
$$ M(a,b,\varepsilon,E,E')=\left[\left( 2 (2 a+b)+\frac1{\kappa_1(\partial E)}+ \delta^H(E,E')\right)^2 -\left(2 (a+b)+\frac1{\kappa_1(\partial E)}\right)^2+(2 b)^2\right]^{\frac 12}+\varepsilon.$$
\end{enumerate}
 Then, there exist a polygon $Q$ and a conformal minimal immersion $Y: \overline{\intc Q} \rightarrow \r^3$, with $Y(0)=0$, and verifying:
\begin{enumerate}[\rm ({d}.1)]
\item $\overline{\intc P^\varepsilon} \subset \intc Q \subset \overline{\intc Q} \subset \intc P$;
\item $\frac 1{\varepsilon}<\dist_{(\overline{\intc Q},\met Y)}(z,P^\varepsilon)$, $\forall z\in Q$;
\item $Y(Q) \subset E' - E'_{-b}$;
\item $Y(\intc Q - \intc P^\varepsilon) \subset \r^3 - E_{-2 (a+b)}$;

\item $\| Y(z)-X(z)\| <\varepsilon$, $\forall z\in \overline{\intc P^\varepsilon}$;

\item  The Hausdorff distance $\delta^H(\Gamma, Y(\overline{\intc Q - \intc P^\varepsilon}))<c+M(a,b,\varepsilon,E,E')$;
\item The cycle  $Y(P)$ is homologous to $\Gamma$  in the open tube 
$T\left( \Gamma,c+M(a,b,\varepsilon,E,E')\right)$
\end{enumerate} 
\end{lema}

%%%%%%%%%%%%%%%%%%%%%%%%%Lema completa %%%%%%%%%%%%%%%%%%%%%%%%%%%%%%%%%%%%%%%%%%%%%%%

\subsection{Proof of Lemma \ref{lem:complete}}

Consider $P$, the polygon given in the statement of the lemma. As usual in constructions of this kind, our first step will consist of describing a labyrinth on $\intc P$, which depends on $P$ and a positive integer $N$. 

Let $\ell$  be the number of sides of $P$. From now on, $N$ will be a positive multiple of $\ell$.

\begin{remark}
Throughout the proof of the lemma a set of real positive constants 
 depending on $X$, $P$, $r$, $\epsilon$ and $s$ will appear. These constants will be represented by the symbol $\cte$. 
Notice that the choice of these constants does not depend 
on the integer $N$.
\end{remark}

Let $\zeta_0>0$ small enough so that $P^{\zeta_0}$ is well defined and 
$\overline{\intc(P^\epsilon)} \subset \intc(P^{\zeta_0})$. From now on, 
we will only consider $N \in \n$ such that $2/N < \zeta_0$.
Let $v_1, \ldots ,v_{2N}$ be a set of points in the polygon $P$ (containing 
the vertices of $P$) that divide each side of $P$ into $\frac{2N}{\ell}$ 
equal parts. We can transfer this partition to the polygon $P^{2/N}$: 
$v_1^\prime, \ldots, v_{2N}^\prime$. We define the following sets:
\begin{itemize}
\item $L_i=$ the segment that joins $v_i$ and $v_i^\prime$, $i=1, \ldots, 2 N$;
\item $\mathcal{P}_i=P^{i/N^3}, \; i=0, \ldots ,2N^2$;
\item $\mathcal{E}= \bigcup_{i=0}^{N^2-1} \overline{\intc(\mathcal{P}_{2i}) 
- \intc(\mathcal{P}_{2i+1})}$ and $\widetilde{\mathcal{E}}= 
\bigcup_{i=1}^{N^2} \overline{\intc(\mathcal{P}_{2i-1}) - \intc(\mathcal{
P}_{2i})}$;
\item $\mathcal R= \bigcup_{i=0}^{2N^2} \mathcal{P}_i$;
\item $\mathcal{B}= \bigcup_{i=1}^N L_{2i}$ and $\widetilde{\mathcal{B}}= 
\bigcup_{i=0}^{N-1} L_{2i+1}$;
\item $L=\mathcal{B} \cap \mathcal{E}$, 
$\widetilde{L}=\widetilde{\mathcal{B}} \cap \widetilde{\mathcal{E}}$, and 
$H=\mathcal R \cup L \cup \widetilde{L}$;
\item $\Xi_N= \{ z \in \intc(\mathcal{P}_0) - \intc(\mathcal{P}_{2N^2}) : 
\dist_{ds_0,\c}(z,H) \geq \frac{1}{4N^3}\}$, where $ds_0$ is the Euclidean 
metric on $\c$.
\end{itemize}
We define $\omega_i$ as the union of the segment $L_i$ and those connected 
components of $\Xi_N$ that have nonempty intersection with $L_i$ for $i=1, 
\ldots,2 N$. Finally, we label $\varpi_i= \{ z \in \c \: : \: 
\dist_{ds_0,\c}(z,\omega_i)< \delta(N) \}$, where   $\delta(N)>0$ is chosen 
in such a way that the sets $\overline{\varpi_i}$ ($i=1, \ldots ,2 N$) are 
pairwise disjoint. 

The shape of the labyrinth formed by the sets $\omega_i$, guarantee the 
following claims if $N$ is large enough:

\begin{claim} The Euclidean 
diameter of $\varpi_i$ is less than $\frac{\cte}{N}$.
\end{claim}

\begin{claim}If $\lambda^2 
\left< \cdot , \cdot \right>$ is a conformal metric on $\overline{\intc P}$ 
and verifies
\[\lambda \geq\begin{cases}
 c & \text{in } \intc P,\\
 c\; N^4 &\text{in } \Xi_N,
\end{cases}\] 
for $c \in \r^+$,
and if $\alpha$ is a curve in $\overline{\intc P}$ connecting 
$P^\epsilon$ and $P$, then
\[\longui (\alpha,\lambda \left<\cdot,\cdot\right>)>\frac{c \cdot \cte \cdot N}{2}.\]
\end{claim}

\noindent Claim  3.2 is a consequence of the fact that a curve $\alpha$, that does 
not go through the connected components of $\Xi_N$, must have a large 
Euclidean length.

Provided that $N$ is large enough, we can assume that $X(\overline{\varpi_i})$ has a sufficiently small diameter in $\r^3$ so that 
\begin{equation}\label{eq:pitonisa}
X(\overline{\varpi_i}) \subset \b(p_i,r), \quad \mbox{ where $p_i \in \Gamma$, $i=1, \ldots 2 N$}.
\end{equation}
For each point $p \in \b(p_i,r)-\{p_i\}$ we define: $$\nor_i(p)=\frac{p-p_i}{\|p-p_i\|}.$$
\begin{figure}
	\begin{center}
		\includegraphics[width=\textwidth]{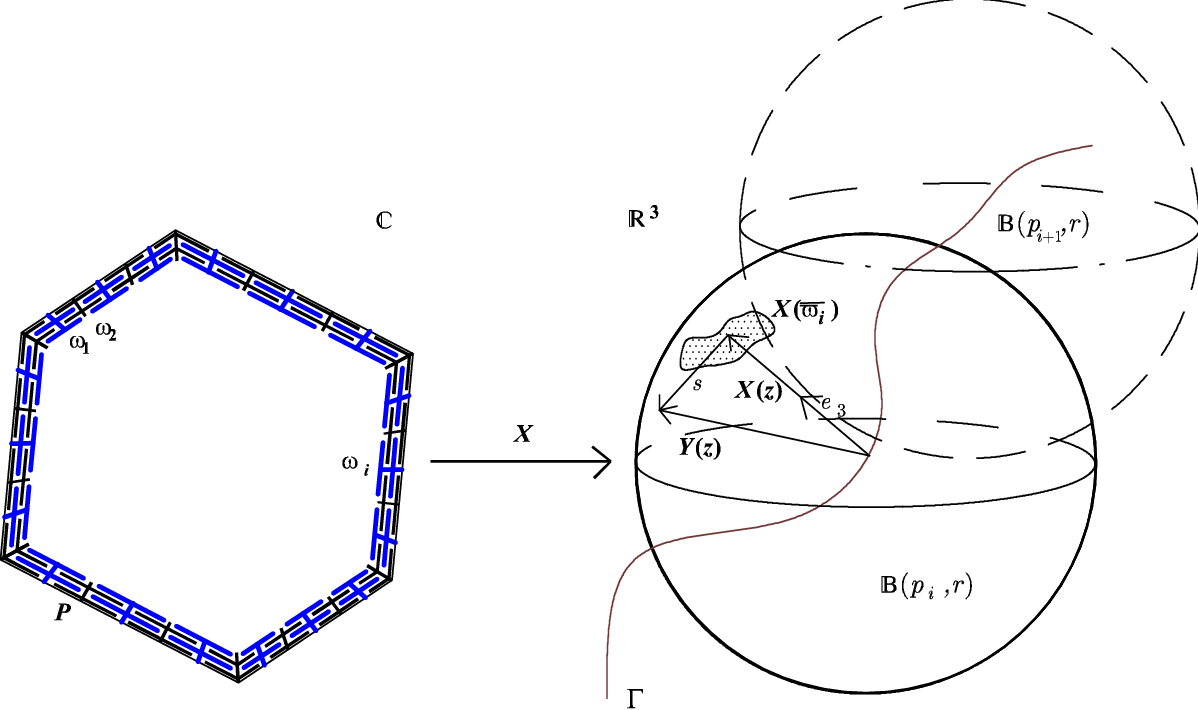}
	\end{center}
	\caption{The set $\varpi_i$ and its image in $\rth.$}
	\label{fig:apli}
\end{figure}

We pursue the construction of a finite sequence of minimal immersions (with boundary); $F_0=X,$ $F_1, \ldots ,F_{2 N},$ satisfying:
\begin{enumerate}[({A}.1{$_i$})]
\item $\|\phi^i(z)-\phi^{i-1}(z)\| \leq \frac{1}{N^2}$, $\forall z \in \overline{\intc P}- \varpi_i;$
\item $\| \phi^i(z)\| \geq N^{7/2},$ $\forall z \in \omega_i;$
\item $\| \phi^i(z) \| \geq \frac{\cte}{\sqrt{N}},$ $\forall z \in \varpi_i;$
\item $\dist_{(\s^2, ds_0)}(G_i(z),G_{i-1}(z)) <\frac 1{N^2},$ $\forall z \in \overline{\intc P} - \varpi_i;$
\item There exists an orthonormal basis of $\r^3$, $S_i=\{e_1^i,e_2^i,e_3^i\}$, so that:
\begin{enumerate}[({A.5}.1{$_i$})]
\item If $z \in \varpi_i$ and $\|X(z)-p_i \|\geq \frac 1{\sqrt{N}},$ then $\| \left(X(z)-p_i\right)_{(*,S_i)} \| < \frac{\cte}{\sqrt{N}};$
\item $\left( F_i(z)\right)_{(3,S_i)}=\left( F_{i-1}(z)\right)_{(3,S_i)}$, $\forall z \in \overline{\intc P};$
\end{enumerate}
\item $\|F_i(z)-F_{i-1}(z) \| \leq \frac{\cte}{N^2},$ $\forall z \in \overline{\intc P} - \varpi_i.$
\item $F_i(P^\epsilon)$ and $\Gamma$ are homologous in $T(\Gamma,R)$
\end{enumerate}

In order to obtain the sequence $F_0, \ldots,F_{2 N}$ we follow an inductive method. Assume
we have constructed $F_0, \ldots, F_{i-1}$ verifying Properties (A.1$_i$), $\ldots$, (A.7$_i$). We define $F_i$ as follows.

First, observe that (for a large enough N) one has:
\begin{enumerate}[({B.}1)]
\item There are positive constants so that $\cte_1 \leq \| \phi^{i-1}(z) \| \leq \cte_2$, for all $z \in \intc P - \cup_{k=1}^{i-1} \varpi_k;$\newline To obtain this property, it suffices to apply (A2$_j$) for $j=1, \ldots, i-1$.
\item The diameter in $\r^3$ of $F_{i-1}(\varpi_i)$ is less than $\frac 1{\sqrt{N}}$. \newline This is a consequence of (B.1), Claim 3.1, and (\ref{eq:metric}).
\item The diameter in $\s^2$ of $G_{i-1}(\varpi_i)$ is less than $\frac 1{\sqrt{N}}$. In particular, the set $G_{i-1}(\varpi_i)$ can be included in a cone $\hbox{Cone}\left(g , \frac{1}{\sqrt{N}}\right)$, for a suitable $g \in G_{i-1}(\varpi_i)$. \newline From Claim 3.1, the diameter of $G_0(\varpi_i)$ is bounded. Then (B.3) holds after successive applications of (A.5$_j$), $j=1, \ldots,i-1.$
\item There exists an orthogonal frame $S_i=\{e_1, e_2, e_3 \}$ in $\r^3$, where: 
\begin{enumerate}[({B.4.}1)] 
\item If $z \in \varpi_i$ and $\|X(z)-p_i \|\geq \frac 1{\sqrt{N}},$ then $\angle \left(e_3, X(z)-p_i\right) \leq \tfrac{\cte}{\sqrt{N}};$
\item $\angle (\pm e_3,G_{i-1}(z)) \geq \tfrac{\cte}{\sqrt{N}}$ for all $z \in \varpi_i$.
\end{enumerate}
The proof of (B.4) is slightly more complicated. Let  $C=\hbox{Cone}\left(g,\frac{2}{\sqrt N}\right)$  where $g$ is given by Property (B.3). To obtain (B.4.2) it suffices to take $e_3$ in  $\s^2 - H$, where
$H=C \bigcup (-C).$
On the other hand, in order to verify (B.4.1), the vector $e_3$ must be chosen as follows:
{\samepage
\begin{itemize}
\item If $(\s^2\setminus H)\cap \nor_i (X(\varpi_i)-\{p_i\}))\not=\emptyset$, then we take $e_3\in(\s^2\setminus H)\cap \nor (X(\varpi_i)-\{p_i\})$; 
\item If $(\s^2\setminus H)\cap \nor_i (X(\varpi_i)-\{p_i\})=\emptyset$, then we take $e_3\in \s^2- H$ satisfying $\angle(e_3,q')<\frac{2}{\sqrt N}$ for some $q'\in \nor_i (X(\varpi_i)-\{p_i\})$.
\end{itemize}
}
It is straightforward to check that this choice of $e_3$ guarantees (B.4).
\end{enumerate}

At this point we are able to construct the $i$-th element of our sequence $F_i$. Let $(g^{i-1}, \phi_3^{i-1})$ be the Weierstrass data of $F_{i-1}$ in the orthonormal basis $S_i$. Applying Runge's theorem, we can construct a family of holomorphic functions, $h_\alpha: \c \rightarrow \c^* $ verifying:
\begin{enumerate}[(a)]
\item $|h_\alpha(z)-1| < \frac 1\alpha$, for all $z \in \overline{\intc P} - \varpi_i;$
\item $|h_\alpha(z)- \alpha| < \frac 1\alpha,$ for all $z \in \omega_i;$
\end{enumerate}
where $\alpha \in \r_+.$ It is important to note that the family $\{h_\alpha|_{\overline{\intc P}- \varpi_i}\; | \; \alpha \in \r_+ \}$ is continuous in the parameter $\alpha.$

Using $h_\alpha$ as a López-Ros function, we define this new Weierstrass data:
\begin{eqnarray*} 
g^{(i,\alpha)} & = & \frac{g^{i-1}}{h_\alpha}, \qquad \phi_3^i \; = \; \phi_3^{i-1},
\end{eqnarray*}
and the associated conformal minimal immersion:
$$F_{(i,\alpha)}(z)= \frac 12 \re \int^z \left( \frac{1}{g^{(i,\alpha)}}-g^{(i,\alpha)} , {\rm i} \left( \frac{1}{g^{(i,\alpha)}}+g^{(i,\alpha)} \right), 2 \right)\phi_3^i .$$
Now, we have to check that there exists a real $\alpha_0>0$ such that $F_i \df F_{(i,\alpha_0)}$ satisfies Properties (A.1$_i$), $\ldots$, (A.7$_i$). 

Since $h_\alpha \to 1$, uniformly on $\overline{\intc P}- \varpi_i$, and $h_\alpha \to \infty$, uniformly on $\omega_i$, as $\alpha \to \infty$, then {\bf (A.1$_i$)}, {\bf (A.2$_i$)}, {\bf (A.3$_i$)} and {\bf(A.6$_i$)} trivially hold for any $\alpha>\alpha_0$, provided $\alpha_0$ is large enough in terms of $N$.

As we mentioned before, the family $\{h_\alpha|_{\overline{\intc P}- \varpi_i}, \; \alpha \in \r_+\}$ depends continuously on $\alpha$, and $h_\alpha \to 1$ as $\alpha \to +\infty$, uniformly on $\overline{\intc P}- \varpi_i$. As $P^\epsilon$ is contained in $\overline{\intc P}- \varpi_i$, then we can see $F_{(i,\alpha)}(P^\epsilon)$, $\alpha \in [\alpha_0,+\infty ]$, as a continuous deformation between $F_{(i,\alpha_0)}(P^\epsilon)$ and $F_{i-1}(P^\epsilon)$. Furthermore, if $N$ is large enough, Property (A.6$_i$) implies that $F_{(i,\alpha)}(P^\epsilon)$ lies in the interior of $T(\Gamma, R)$, for all $\alpha \in [\alpha_0,+\infty].$ So, Property {\bf (A.7$_i$)} is a consequence of these facts and Property (A.7$_{i-1}$).

In order to check {\bf (A.3$_i$)} we have to use (B.4.2). This property gives us:
$$\frac{\sin \left(\frac{\cte}{\sqrt{N}}\right)}{1+\cos\left(\frac{\cte}{\sqrt{N}}\right)} \leq |g^{i-1}| \leq \frac{\sin \left(\frac{\cte}{\sqrt{N}}\right)}{1-\cos\left(\frac{\cte}{\sqrt{N}}\right)} \qquad \hbox{in } \varpi_i,$$
and so, taking (B.1) into account one has (if $N$ is large enough):
$$\| \phi^i \| \geq | \phi^i_3|=| \phi^{i-1}_3| \geq \sqrt2\| \phi^{i-1} \| \frac{|g^{i-1}|}{1+|g^{i-1}|^2}\geq \cte \cdot \sin\left(\tfrac{\cte}{\sqrt{N}}\right) \geq \tfrac \cte{\sqrt{N}} \qquad \hbox{in } \varpi_i.$$

Using (B.4.1), we get {\bf (A.5.1$_i$)}. Finally, to obtain {\bf (A.5.2$_i$)}, we use that $\phi^{i-1}_3=\phi_3^i$ in the frame $S_i$.
\vskip 1cm

Hence, we have constructed the immersions $F_0,F_1, \ldots ,F_{2N}$ verifying claims (A.1$_i$),$\ldots$,(A.7$_i$) for $i=1, \ldots , 2N$. The following proposition stands all the properties of $F_{2 N}$ we will need.
\begin{proposicion} \label{prop:F2N}
If $N$ is sufficiently large, then we have:
\begin{enumerate}[\rm (I)]
\item $\dist_{( \overline{\intc P}, F_{2 N})} (P, P^\epsilon)>2 s ;$
\item $\|F_{2 N}(z)-X(z)\| < \cte/N$, for all $\displaystyle z \in \overline{\intc P} - \bigcup_{i=1}^{2 N} \varpi_i;$
\item There exists a polygon $\widetilde P$ in $\c$ such that:
\begin{enumerate}[\rm {(III.} a{)}]
\item $\overline{\intc P^\epsilon} \subset \intc{\widetilde{P}} \subset \overline{\intc{\widetilde{P}}} \subset \intc{P}$;
\item  $s<\dist_{(\overline{\intc P},\metri{F_{2N}})}(z,P^\epsilon)<2s$, $\forall z \in \widetilde{P}$;
\item $F_{2N}(\overline{\intc P- \intc P^\epsilon}) \subset T(\Gamma, R)$;
\end{enumerate}
\item $F_{2 N}(P^\epsilon)$ and $\Gamma$ are homologous in $T(\Gamma, R)$.
\end{enumerate}
\end{proposicion}
\begin{proof}
Item (I) is an standard consequence of Claim 3.2, making use of (A.1$_i$), (A.2$_i$) and (A.3$_i$), $i=1, \ldots, 2N.$ Similarly, a successive application of (A.6$_i$) implies item (II).

The demonstration of item (III) is a bit more delicate. First we need to construct the polygon $\widetilde P$ in (III). To do this, we consider the set ${\cal A}= \{ z \in \intc P - \intc P^\epsilon \; : \; s<\dist_{(\overline{\intc P},\metri{F_{2N}})}(z,P^\epsilon)<2s\}$. Note that $\cal A$ is nonempty and that $P$ and $P^\epsilon$ are in different connected components of $\c - {\cal A}.$ Then, the existence of the  polygon $\widetilde P$ satisfying items (III.a) and (III.b) is obvious.

Now, we are going to check that $F_{2N}$ verifies item (III.c). Consider $\eta \in \overline{\intc \widetilde P - \intc P^\epsilon}.$ We will assume that $F_{2 N}(\eta) \not\in \Gamma$, otherwise we have nothing to prove. Hence, we distinguish two cases:

\noindent {\sc \bf Case 1.} The point $\eta$ belongs to $\displaystyle \intc P - \bigcup_{i=1}^{2 N} \varpi_i.$ 

In this case we know that $\|F_{2 N}(\eta)-X(\eta) \| \leq \frac{\cte}{N}.$ Since $X(\eta)$ lies in the interior of $T(\Gamma,r)$, then we can choose $N$  large enough so that $F_{2N}(\eta) \in T(\Gamma,r) \subset T(\Gamma, R).$
\vskip .5cm

\noindent {\sc \bf Case 2.} There exists $i \in \{1, \ldots, 2N\}$ such that $\eta \in \varpi_i.$

Consider now a curve $\beta:[0,1] \rightarrow \intc P$ so that $\beta(0) \in P^\epsilon$, $\beta(1)=\eta$ and $\longui(\beta, F_{2 N}) \leq 2 s.$ Note that the existence of such a curve is guaranteed by (III.b). 
Let us define $\overline{t}= \mbox{Supremum} \{ t \in [0,1] \; : \; \beta(t) \in \partial \varpi_i\}$ and $\overline \eta = \beta(\overline t).$ It is important to note that $\overline t$ exists because $\varpi_i \subset \intc P - \intc P^\epsilon.$

For our purposes, we need to prove first the following inequality:
\begin{equation}\label{eq:lembra}
\| F_i(\overline \eta)-F_i(\eta)\| \leq \frac{\cte}{N}+ 2 s.
\end{equation}
Indeed, by using Properties (A.6$_k$), for $k=1, \ldots , 2N$, we obtain:
\begin{multline*}
\| F_i(\overline \eta)-F_i(\eta) \| \leq \| F_i(\overline \eta)-F_{2 N}(\overline \eta) \|+ 
\| F_{2 N}(\overline \eta)-F_{2 N}(\eta)\|+\|F_{2 N}(\eta)-F_i(\eta)\| \leq \\
\frac{\cte}{N}+\longui(\beta,F_{2 N})+\frac{\cte}{N} \leq \frac{\cte}{N}+ 2 s.
\end{multline*}
Let $p_i$ be the point given by condition \eqref{eq:pitonisa}. Then we have:
\begin{equation} \label{eq:correa}
\|F_{2 N}(\eta)-p_i\|\leq \|F_i(\eta)-p_i\|+ \frac{\cte}{N}.
\end{equation}
We again distinguish two cases:

\noindent {\sc \bf Case 2.1.} If  $\| X(\eta)-p_i\| \leq 1/\sqrt{N}$, then one has:
\begin{multline*}
\|F_i(\eta)-p_i\|\leq \| F_i(\eta)-F_i(\overline \eta)\|+\|F_i(\overline \eta)-F_{i-1}(\overline \eta)\|+
\|F_{i-1}(\overline \eta)-F_{i-1}(\eta) \|+\|F_{i-1}(\eta)-X(\eta)\|+ \|X(\eta)-p_i \| \leq \\
\frac{\cte}N+2 s+\frac{\cte}{N^2}+\frac{\cte}{\sqrt{N}}+\frac{\cte}N+\frac {1}{\sqrt{N}} \leq R
\end{multline*}
where the last inequality occurs if $N$ is sufficiently large.

\noindent {\sc \bf Case 2.2.} If  $\| X(\eta)-p_i\| > 1/\sqrt{N}$, then we use (A.5.2$_i$) to get a bound for the third coordinate of $F_i(\eta)-p_i$ in the orthonormal basis $S_i$. We proceed as follows: 
\begin{multline} \label{eq:tres}
|(F_i(\eta)-p_i)_{(3,S_i)}|=|(F_{i-1}(\eta)-p_i)_{(3,S_i)}|\leq |(F_{i-1}(\eta)-X(\eta))_{(3,S_i)}|+|(X(\eta)-p_i)_{(3,S_i)}| \leq 
\frac{\cte}{N}+r.
\end{multline}
On the other hand, we can apply Property (A.5.1$_i$) to find an upper bound for the first two coordinates of $F_i(\eta)-p_i,$
\begin{multline}\label{eq:dos}
\| (F_i(\eta)-p_i)_{(*,S_i)} \| \leq \|(F_i(\eta)-F_i(\overline \eta))_{(*,S_i)} \| +\|(F_i(\overline \eta)-F_{i-1}(\overline \eta))_{(*,S_i)}\|+ \\
\|(F_{i-1}(\overline \eta)-F_{i}( \eta))_{(*,S_i)}\|+\|(F_{i-1}(\eta)-X( \eta))_{(*,S_i)}\|+\|(X(\eta)-p_i)_{(*,S_i)}\| \leq \\
\frac{\cte}{N}+2 s+ \frac{\cte}{N^2}+\frac{1}{\sqrt{N}}+\frac{\cte}{\sqrt{N}} \leq 2 s+\frac{\cte}{\sqrt{N}}.
\end{multline}

By Pythagoras' theorem and taking into account (\ref{eq:tres}) and (\ref{eq:dos}), we easily infer that:
$$\| F_i(\eta)-p_i\| < \sqrt{\left(2 s+ \frac{\cte}{\sqrt{N}}\right)^2+ \left(r+\frac{\cte}{N}\right)^2}.$$
Using this upper bound, inequality (\ref{eq:correa}) becomes:
$$\|F_{2 N}(\eta)-p_i\| \leq \sqrt{\left(2 s+ \frac{\cte}{\sqrt{N}}\right)^2+ \left(r+\frac{\cte}{N}\right)^2}+ \frac{\cte}{N}.$$
So, for a large enough $N$, it is clear that $F_{2 N}(\eta) \in T(\Gamma,R)$ where $R=\sqrt{r^2+(2 s)^2}+\epsilon$, for all $\eta \in \overline{\intc \widetilde P- \intc P^\epsilon}.$ This completes the proof of item (III.c).

Finally, item (IV) in the proposition is a direct consequence of property (A.7$_{2 N}$).
\end{proof}

Obviously, the immersion $F_{2N}:\intc \widetilde P \rightarrow \r^3$ is the immersion $\widetilde X$ we are looking for. Proposition \ref{prop:F2N} says to us that $\widetilde X$ so defined verifies items (a.1), $\ldots$, (a.6) in the 
lemma.

%%%%%%%%%%%%%%%%%%%%%%%%%MEEKS' TRICK %%%%%%%%%%%%%%%%%%%%%%%%%%%%%%%%%%%%%%%%%%%%%%%%

\subsection{Proof of Theorem \ref{th:meeks}}

Let $c_0>0$, $r_1>0$ and $\rho_1>0$ to be specified later, and define $$r_n=\sqrt{r_{n-1}^2+\left(\frac{2 c_0}{n}\right)^2}+\frac{c_0}{n^2},$$  and $\rho_n=\rho_1+\sum_{i=2}^n c_0/i,$ $n \geq 2$. The constants $r_1$ and $c_0$ are choosen in such a way that:
\begin{eqnarray}
\lim_{n \to \infty} r_n & < & \frac{\mu}4, \label{eq:bonete-1} \\
\sum_{n=1}^\infty \frac{c_0}{n^2} & < & \frac{\mu}{4}. \label{eq:bonete-2}
\end{eqnarray}
Consider also a regular curve $\Gamma$ in $\r^3$, so that $\delta^H(X(P),\Gamma) <r_1$ and both curves are homologous in $T(\Gamma,r_1).$
Our strategy consists of using Lemma \ref{lem:nadi} to define  a sequence:
$$\chi_n=(X_n:\overline{\intc P_n} \rightarrow \r^3, P_n, \varepsilon_n,\xi_n),$$
where $X_n$ is a conformal minimal immersion,  $P_n$ is a polygon,  $\{ \varepsilon_n \}$, $\{ \xi_n \}$ are decreasing sequences of non vanishing terms satisfying $\varepsilon_n,\xi_n<c_0/n^2$, and:
\begin{enumerate}[$(A_{n})$]
\item  $\overline{\intc P_{n-1}^{\xi_{n-1}}}\subset \intc P_{n-1}^{\varepsilon_n}\subset \overline{\intc P_{n-1}^{\varepsilon_n}}\subset \intc P_n^{\xi_n} \subset \overline{ \intc P_n^{\xi_n}} \subset \intc P_n\subset  \overline{\intc P_n}\subset \intc P_{n-1},$
\item  $\rho_{n}<\dist_{(X_{n},\overline{\intc P_n^{\xi_n}})}(0,P_n^{\xi_n}),$
\item  $\| X_n-X_{n-1}\| < \varepsilon_n$ in $\intc P_{n-1}^{\varepsilon_{n}},$
\item  ${\displaystyle \lambda_{X_n} \geq \alpha_n \lambda_{X_{n-1}}}$ in $\intc P_{n-1}^{\xi_{n-1}},$
where $\{\alpha_i\}_{i \in \n}$ is a sequence of real numbers such that $0<\alpha_i<1$ and $\{ \prod ^n_{i=1} \alpha_i\}_n$ converges to $1/2$,
\item $X_{n}(P_{n-1}^{\varepsilon_n})$ and $\Gamma$ are homologous in $T(\Gamma,r_n),$
\item $X_n\left(\overline{\intc P_n - \intc P_{n-1}^{\varepsilon_n} }\right) \subset T \left( \Gamma, r_n \right).$

\end{enumerate}
The choice of the first element of the sequence is not difficult. First, we take $X_1=X.$ Let $P_1$ be a polygon parallel to $P$ and satisfying $P_1^{\xi_1} =P$. The constant and the polygon must be chosen in such a way that  $\delta^H(X_1 \left(P_1 - P_1^{\xi_1}\right),\Gamma) < r_1.$  Finally, we choose $\rho_1$ and $\varepsilon_1$ satisfying: $$\rho_1<\dist_{(X_1,\overline{\intc P_1})}(0,P_1) \quad \mbox{and} \quad \varepsilon_1<\min \{c_0,r_1\}.$$
Suppose that we have defined $\chi_1, \ldots, \chi_n$. Then, we will construct the $(n+1)$-th term in the following way. 

%Choose $\xi_{n+1}<1/(n+1)^2$ verifying  (B$_{n+1}$) and $\overline{\intc(P_{n-1}^{\xi_{n}})} \subset \intc(P_{n}^{\xi_{n+1}})$, (the choice of $\xi_{n+1}$ is possible since $\chi_n$ satisfies (A$_n$) and  (G$_n$)).

Take a sequence $\{\widehat{\varepsilon}_m \}\searrow 0$, with $\widehat{\varepsilon}_m<\frac{c_0}{(n+1)^2}$, $\forall m$. For each $m$, we consider $\widetilde{P}_m$ and $Y_m:  \intc \widetilde P_m\rightarrow \r^3$  given by Lemma \ref{lem:complete}, for the data:
$$X=X_n, \; P=P_n, \; r=r_n, \; s=\frac{c_0}{n+1}, \; \nu=\widehat{\varepsilon}_m.$$
If $m$ is large enough, items (a.1) and (a.2) in Lemma \ref{lem:complete} tell us that $\overline{\intc P_n^{\xi_{n}}} \subset  \intc\widetilde P_m$ and the sequence $\{Y_m\}$ converges to $X_n$ uniformly in $\overline{\intc P_n^{\xi_{n}}}$. In particular, $\{\lambda_{Y_m}\}$ converges uniformly to $\lambda_{X_n}$  in $\overline{\intc P_n^{\xi_{n}}}$. Therefore there is a $m_0 \in \n$ such that:
\begin{eqnarray} 
\overline{\intc P_n^{\xi_{n}}} \subset & \intc P_n^{\hat \varepsilon_{m_0}}& \subset  \intc \widetilde P_{m_0}, \label{clavao1-p}\\
\rho_{n} & < & \dist_{(Y_{m_0},\overline{\intc P_n^{\xi_{n}}})}(0, P_n^{\xi_{n}}), \label{clavao-p} \\
\lambda_{Y_{m_0}}& \geq & \alpha_{n+1} \lambda_{X_n} \qquad \hbox{in } \intc P_n^{\xi_{n}}.\label{lambdas-p}
\end{eqnarray}
We define $X_{n+1}=Y_{m_0}$, $P_{n+1}=\widetilde{P}_{m_0}$, and $\varepsilon_{n+1}=\widehat{\varepsilon}_{m_0}$. From (\ref{clavao1-p}), (\ref{clavao-p}) and Item (a.3) in Lemma  \ref{lem:complete}, it is not hard to see that $\rho_{n+1}< \dist_{(X_{n+1},\overline{\intc P_{n+1}})}(0,P_{n+1})$. Finally, take $\xi_{n+1}$ small enough such that (A$_{n+1}$) and (B$_{n+1}$) hold. The remaining properties directly follow from (\ref{clavao1-p}), (\ref{lambdas-p}) and the aforementioned lemma. This concludes the construction of the  sequence $\{\chi_n\}_{n \in \n}$.

Now, we define 
$$\Sigma=\bigcup_{n=1}^\infty \intc(P_n^{\xi_{n}})\left(= \bigcup_{n=1}^\infty \intc(P_n^{\varepsilon_{n+1}}) \right).$$
$\Sigma$ is a simply connected domain in $U$. 
Properties (C$_n$) and the fact that $\varepsilon_n<c_0/n^2$ give us that the sequence of minimal immersions $\{ X_n \}$ is a Cauchy sequence, uniformly on compact sets of $\Sigma$, and so $\{ X_n \}$ converges.

Let $\widehat X:\Sigma \rightarrow \r^3$ be the limit of $\{ X_n \}$. $\widehat X$ has the following properties:
\begin{itemize}

\item $\widehat X$ is an immersion. Indeed, for any $z \in \Sigma$ there exists  $n \in \n$ such that $z \in \intc P_n^{\xi_{n}}$. From Properties (F$_i$), $i=n+1, \ldots,k $ we get:
$$ \lambda_{X_k}(z) \geq \alpha_k \lambda_{X_{k-1}}(z) \geq \ldots \geq \alpha_k \ldots \alpha_{n+1} \lambda_{X_n}(z) \geq \alpha_k \ldots \alpha_1 \lambda_{X_n}(z), \; \forall k>n.$$
Taking limit as $k \to \infty$, we deduce:
\begin{equation} \label{immersion}
\lambda_{\widehat X}(z) \geq \frac{1}{2} \; \lambda_{X_n}(z) > 0,
\end{equation}
and so $\widehat X$ is an immersion.
\item $\widehat X$ is minimal and conformal.
\item $\Sigma$ is complete  with the metric induced by $\widehat X$. Indeed, if $n$ is large enough, and  taking (\ref{immersion}) and (A$_n$) into account, one has: 
$$\dist_{(\widehat X,\overline{\intc P_n^{\xi_{n}}})}(0,P_n^{\xi_{n}})>\frac12 \dist_{(X_n,\overline{\intc P_n^{\xi_{n}}})}(0, P_n^{\xi_{n}}) > \frac12\rho_n.$$
The completeness is due to the fact that  $\{\rho_n\}_{n \in \n}$ diverges.
\item  $ \widehat X(\Sigma - \intc P) \subset T(X(P), \mu).$ Pick a point $z \in \Sigma - \intc P$. Then we know that $z$ belongs to $\overline{\intc P_n - \intc P_{n-1}^{\varepsilon_n}}$, for some $n \in \n,$ and so Property (F$_n$) implies $\dist_{\r^3}(X_n(z), \Gamma) <r_n.$ Therefore, one has:
$$
\dist_{\r^3}(\widehat X(z), \Gamma) \leq \|\widehat X(z)-X_n(z)\|+r_n \leq 
\left(\sum_{k=n}^\infty \|X_{k+1}(z)-X_k(z) \| \right)+r_n,$$
at this point we use Properties (C$_k$), for $k \geq n$, and then one obtains:
$$\dist_{\r^3}(\widehat X(z), \Gamma) \leq \sum_{k=n}^\infty \varepsilon_k + r_n < \frac{\mu}2.$$
Thus we have that $\widehat X(\Sigma - \intc P) \subset T (\Gamma,\mu/4).$ Moreover, our choice of $\Gamma$ implies that $\widehat X(\Sigma - \intc P) \subset T (X(P),\mu/2+r_1) \subset  T (X(P),\mu).$
\item $\widehat X(P)$ is homologous to $X(P)$ in $T(X(P), \mu).$ Indeed, pick a natural $n \in \n,$ then from Properties C$_k$, $k=1, \ldots, n$ we have that $X_n(P)\subset T(\Gamma,r_n).$ It is clear that $X_n(P)$ is homologous to $X_n(P_{n-1}^{\varepsilon_n})$ in $T(\Gamma,r_n)$, so using Property (E$_n$) and taking into account our choice of $\Gamma$ we conclude that $X_n(P)$ and $X(P)$ are homologous in $T(\Gamma,r_n) \subset T(X(P),\mu).$ Since the curves $X_n(P)$ converge uniformly to $\widehat X(P)$ and all these curves have the same homological type as $X(P)$, then  it is clear that $\widehat X(P)$ is also homologous to $X(P)$ in the tube $T(X(P),\mu)$. 
\end{itemize}This completes the proof of Theorem \ref{th:meeks}
%%%%%%%%%%%%%%%%%%%%%%%%%LEMA DE NADIRASHVILI %%%%%%%%%%%%%%%%%%%%%%%%%%%%%%%%%%%%%%%%%%%%%%%%%%%%

\subsection{Proof of Lemma \ref{lem:nadi}}

All the arguments we need to prove this lemma are essentially contained in Theorem \ref{th:meeks} and Lemma \ref{lem:propia}.
First, we apply Theorem \ref{th:meeks} to the immersion $X: \intc P \rightarrow \r^3$, the polygon $P^\varepsilon$ and a constant $\mu>0$, to be determined later. Hence, we obtain a complete minimal immersion  $\widehat X:\Sigma \rightarrow \r^3$ which satisfies:
\begin{enumerate}[{(a.}i)]
\item $\overline{\intc P^\varepsilon} \subset \Sigma \subset \overline \Sigma \subset \intc P,$
\item $\|X(z)-\widehat X(z)\| < \mu, $ $\forall z \in \intc P^\varepsilon,$
\item $\widehat X(\Sigma - \intc P^\varepsilon ) \subset T(X(P^\varepsilon),\mu),$
\item $\widehat X(P^\varepsilon)$ is homologous to $X(P^\varepsilon)$ in $T(X(P^\varepsilon), \mu).$ 
\end{enumerate}

As the immersion $\widehat X$ is complete, then we can find a new polygon $\widehat P$ satisfying $$\overline{\intc P^\varepsilon} \subset \intc \widehat P \subset \overline{\intc \widehat P} \subset \intc P,$$ and so that the distance \begin{equation}\label{eq:1/e}
\dist_{(\overline{\intc \widehat P},\widehat X)}(z, P^\varepsilon ) >\frac 1{\varepsilon}, \quad \forall z \in \widehat P.
\end{equation}
Taking into account the hypotheses  of the lemma, Property (a.iii) guarantees (for a small enough $\mu$) that $\widehat X(\intc \widehat P-\intc P^\varepsilon)\subset E - E_{-a}$ and $\widehat X(\intc \Sigma-\intc P^\varepsilon)\subset T(\Gamma,c).$ The former inclusion yields that $\delta^H(\widehat X(\intc \Sigma-\intc P^\varepsilon), \Gamma) <c.$ Indeed, we proceed by contradiction. Suppose this is not true, then there exists a point $x \in \Gamma$ so that $\b(x,c) \cap \widehat X(\intc \widehat P-\intc P^\varepsilon)$ is empty. Among other things, this implies that $\widehat X(P^\varepsilon)$ can not be homologous to $\Gamma$ in $T(\Gamma,c)$ which is contrary to Property (a.iv) (recall that $X(P^\varepsilon)$ is homologous to $\Gamma.$)

At this point, we are able to apply Lemma \ref{lem:propia} to the following data:
$$X=\widehat X, \quad P=\widehat P, \quad O=\intc P, \quad E=E_{-a}, \quad E'=E', \quad b_1 < \varepsilon/2 , \quad b_2=b.$$
Thus, we obtain a new polygon $Q$ such that $\intc \widehat P\subset \overline{\intc \widehat P} \subset \intc Q \subset \overline{\intc Q} \subset \intc P$ and a new minimal immersion $Y:\overline{\intc Q}\rightarrow \r^3$ with the following properties:
\begin{enumerate}[({b}.i)]
\item $\|Y(z)- \widehat X(z)\|<b_1$, $\forall z \in \overline{\intc \widehat P},$
\item $Y(Q) \subset E' - E'_{-b},$
\item $Y(\intc Q\setminus \intc \widehat P) \subset \r^3-  E_{-2 b-a},$
\item For any $z\in \intc Q\setminus\intc \widehat  P$, one has: \begin{multline*}\| \widehat X(z)-Y(z) \| <M(a,b,b_1,E,E')= \left[\left( 2 (2 a+b)+\frac1{\kappa_1(\partial E)}+ \delta^H(E,E')\right)^2 -\right. \\ 
\left.-\left(2 (a+b)+\frac1{\kappa_1(\partial E)}\right)^2+(2 b)^2\right]^{\frac 12}+b_1;\end{multline*}
\item $\widehat X(\widehat P)$ and $Y(\widehat P)$ are homologous in $T(\widehat X(\widehat P),M(a,b,b_1,E,E')).$
\end{enumerate}

To finish the proof, it suffices to check that $Y$ is the immersion we are looking for. First, observe that item (d.1) trivially holds.

In order to check (d.2), notice that $Y$ converges to $\widehat X$, uniformly on $\overline{\intc \widehat P}$ as $b_1 \to 0$. So, item (d.2) is an easy consequence of this fact and \eqref{eq:1/e}.

Item  (d.3) directly follows from (b.ii). Moreover, Properties (a.ii), (b.i), and (b.iii) imply item (d.4) (provided $\mu$ and $b_1$ are sufficiently small.)

Similarly, if $\mu,$ $b_1 < \varepsilon/2$, then Item (d.5) is obtained from (a.ii) and (b.i).

Let us show that immersion $Y$ satisfies item (d.6). We know that $$\delta^H \left(\Gamma,\widehat X(Q - \intc P^\epsilon ) \right)< c .$$  Thus, Item (d.6) follows from (b.i) and (b.iv), provided that  $b_1$ is small enough.

Finally, item (d.7) in the lemma can be easily deduced from (a.iv) and (b.v), taking into account that $T(\widehat X(\widehat P),M(a,b,b_1,E,E'))$ and $ T(X(P^\varepsilon), \mu)$ are contained in the tube $T(\Gamma, c+M(a,b,b_1,E,E'))$ and the fact  that $X(P^\varepsilon)$ is homologous to $\Gamma$ in that tube.
This completes the proof of Lemma \ref{lem:nadi}.
%%%%%%%%%%%%%%%%%%%%%%%%%%%%%%%%%%%%%%%%%%%%%%%%%%%%%%%

\section{Proof of the main Theorems}\label{sec:main}
At this point we are able to prove the principal results of our paper. From now on $C$ will represent a bounded, strictly convex  regular domain of $\rth.$ Recall that strictly convex means that the mean and Gaussian curvatures are positive. If we call $(\cal K, \delta^H)$ as the metric space of compact sets of $\partial C$ with the Hausdorff distance, then we will prove that the limit sets of properly immersed minimal disk are dense in $\cal K.$ At the end of this section we will show that our results are sharp in the sense that neither strictly convexity nor regularity can be removed from our assumptions.
\begin{teorema} \label{th:mari}
Let $C$ be a  strictly convex bounded regular domain of space. For any smooth Jordan curve $\Gamma \subset \partial C$ and for any $\epsilon>0$ there exists a complete proper minimal immersion $\psi_{(\Gamma,\epsilon)}:\d \rightarrow C$ so that $\delta^H( \psi_{(\Gamma,\epsilon)}(\partial \d), \Gamma)<\epsilon$.
\end{teorema}

\begin{proof} We know that there is $t_0>0$ so that $C_t$ is well defined for any $t \in ]-t_0, \infty[$. Furthermore, the normal projection $\pro_t:\r^3 - C_t \rightarrow \partial(C_t)$ is a well defined smooth map. 

We fix $\tau >0 $  such that $\tau < \min \{t_0, \epsilon\}$. 
Now, define
$$M(n,c) \df c \left( \frac{1}{n^2}+\sqrt{\frac{29 \, c^2}{(n-1)^8}+\frac{4}{\kappa_1(\partial C) \, (n-1)^4}} \; \right).$$
Choose $c_1>0$ small enough such that $\sum_{n\geq 2}M(n,c_1) <\tau/4$. In particular, $$ \sum_{n = 1}^{\infty} \frac{c_1^2}{n^4} < \sum_{n=2}^{\infty}M(n,c_1) <\tau/4<t_0.$$
Using this constant $c_1$, we construct an expansive sequence $\{E^n \}_{n \in \n}$ of bounded convex regular domains in the following way: $E^n\df C_{-t_n}$, where $t_n=\sum_{k=n}^\infty c_1^2/k^4$, $n\geq 1$.
If we label  $\Gamma_1 \df \pro_{-t_1}(\Gamma)$, it is obvious that $\Gamma_1 \subset \partial E^1$ is a Jordan curve.

We also take a decreasing sequence of positive reals $\{b_n\}_{n \in \n}$, satisfying $b_n<c_1^2/n^4$, $\forall n \in \n$.

Finally, we define a sequence of real numbers $\{ \mu_n \}_{n \in \n}$ in the following way: $\mu_1=\tau/4$ and $\mu_n=\mu_{n-1}+M(n,c_1)$.

The next step consists of using Lemma \ref{lem:nadi} to construct (in a recursive way) a sequence 
$$\chi_n=(X_n:\overline{\intc P_n} \rightarrow \r^3,P_n, \varepsilon_n,\xi_n),$$
where $X_n$ are conformal minimal immersions with $X_n(0)=0$, $P_n$ are polygons,  and $ \{     \varepsilon_n\}$, $\{\xi_n\}$, $\{\sigma_n\}$ are sequences of positive numbers converging to zero, verifying $\varepsilon_k<c_1/k^2$, $\sum_{k=1}^\infty \varepsilon_k <\delta$.
Furthermore, the sequence $X_n:\overline{\intc P_n} \rightarrow \r^3$ must verify the following properties:
\begin{enumerate}[(I$_{n}$)]
\item $\overline{\intc P_{n-1}^{\xi_{n-1}}}\subset \intc P_{n-1}^{\varepsilon_n}\subset \overline{\intc P_{n-1}^{\varepsilon_n}}\subset \intc P_n^{\xi_n} \subset \overline{ \intc P_n^{\xi_n}} \subset \intc P_n\subset  \overline{\intc P_n}\subset \intc P_{n-1}$;

\item $\| X_n(z)-X_{n-1}(z)\| < \varepsilon_n$, $\forall z\in\intc P_{n-1}^{\varepsilon_{n}}$;

\item $\lambda_{X_n}(z) \geq \alpha_n \lambda_{X_{n-1}}(z)$, $\forall z\in \intc P_{n-1}^{\xi_{n-1}}$, where $\{\alpha_i\}_{i \in \n}$ is a sequence of real numbers such that $0<\alpha_i<1$ and $\{ \prod ^n_{i=1} \alpha_i\}_n$ converges to $1/2$;

\item $\displaystyle \frac1{\varepsilon_n} < \dist_{(\overline{\intc P_n^{\xi_n}},X_n)}(P_{n-1}^{\xi_{n-1}},P_n^{\xi_n})$;

\item $X_n(z)\in E^{n}- (E^{n})_{-b_n}$, for all $z\in P_n$;
\item $X_n(z)\in \r^3\setminus (E^{n-1})_{-2(b_{n-1}+b_n)}$, for all $z\in \intc P_n\setminus \intc P_{n-1}^{\varepsilon_n}$;
\item $\displaystyle \delta^H \left(X_n(P_n), \Gamma \right)< \mu_n;$
\item $\displaystyle \delta^H \left(X_n(\intc P_n - \intc P_{n-1}^{\varepsilon_n}),\Gamma \right)< \mu_n;$
\item $X_{n}(P_n)$ and $\Gamma$ are homologous in $T\left(\Gamma, \mu_n \right).$
\end{enumerate}

To define $\chi_1$, we consider $D_1$ a solution of Plateau's problem for the curve $\Gamma_1$. Let $X_1: \overline \d \rightarrow D_1$ be a conformal parametrization of the minimal disk $D_1$. Then we choose a polygon $P_1 \subset \d$ sufficiently close to $\partial \d$ and a constant $\xi_1>0$ so that:
\begin{itemize}
\item  $\delta^H(X_1(\intc P_1 - \intc P_1^{\xi_1}), \Gamma_1)<b_1$ ($< \frac{\tau}{4}$.) In particular, Properties (V$_1$) and (VII$_1$) hold.
\item $X_1(P_1)$ is homologous to $\Gamma$ in the tube $T(\Gamma,\mu_1).$
 \end{itemize}

Suppose that we have $\chi_1, \ldots, \chi_n$. In order to construct $\chi_{n+1}$, we consider the following data:
$$E=E^n,\quad E'=E^{n+1},\quad a=b_{n},\quad c=\mu_n, \quad X=X_n,\quad P=P_n.$$
Property (V$_n$) says to us that $X(P)\subset E\setminus E_{-a}$, and Property (VII$_n$) says that $\delta^H(X(P),\Gamma)< c$. Furthermore, $X(P)$ is homologous to $\Gamma$ in $T(\Gamma,c)$. Then it is straightforward that we can find a small enough positive constant $\varkappa$, such that Lemma \ref{lem:nadi} can be applied to the aforementioned data, and for any $\varepsilon \in ]0,\varkappa[$.

Take a sequence $\{\widehat{\varepsilon}_m \}\searrow 0$, with $\widehat{\varepsilon}_m<\mbox{minimum}\{\frac{c_1}{(n+1)^2},\varkappa , b_{n+1}\}$, $\forall m$. 
For each $m$, we consider $Q_m$ and $Y_m: \overline{\intc Q_m} \rightarrow \r^3$  given by Lemma \ref{lem:nadi}, for the above data and $\varepsilon=b=\widehat{\varepsilon}_m$. It is important to note that the constant $M(a,b,\widehat{\varepsilon}_m,E,E')$ in Lemma \ref{lem:nadi} is less that $M(n+1,c_1).$

If $m$ is large enough, Assertions (b.1) and (b.5) in Lemma \ref{lem:nadi} tell us that $\overline{\intc P_n^{\xi_{n}}} \subset \intc Q_m$ and the sequence $\{Y_m\}$ converges to $X_n$ uniformly in $\overline{ \intc P_n^{\xi_{n}}}$. In particular, $\{\lambda_{Y_m}\}$ converges uniformly to $\lambda_{X_n}$  in $\overline{ \intc P_n^{\xi_{n}}}$. Therefore there is a $m_0 \in \n$ such that:
\begin{eqnarray} 
\overline{ \intc P_n^{\xi_{n}}} \subset &  \intc P_n^{\hat \varepsilon_{m_0}}& \subset \intc Q_{m_0}, \label{clavao1}\\
\lambda_{Y_{m_0}}& \geq & \alpha_{n+1} \lambda_{X_n} \qquad \hbox{in } \intc P_n^{\xi_{n}}.\label{lambdas}
\end{eqnarray}
We define $X_{n+1}=Y_{m_0}$, $P_{n+1}=Q_{m_0}$, and $\varepsilon_{n+1}=\widehat{\varepsilon}_{m_0}$. From (\ref{clavao1}) and Statement (d.2) in Lemma \ref{lem:nadi}, we infer that $\frac{1}{\varepsilon_{n+1}} < \dist_{\left(\overline{ \intc P_{n+1}}, X_{n+1} \right)}(P_n^{\xi_n},P_{n+1})$. Finally, take $\xi_{n+1}$ small enough such that (I$_{n+1}$) and (IV$_{n+1}$) hold.
 
The remainder properties directly follow from (\ref{clavao1}), (\ref{lambdas}) and Lemma \ref{lem:nadi}. 
This concludes the construction of the  sequence $\{ \chi_n \}_{n \in \n}$.
\bigskip

Now, we extract some information from the properties of $\{\chi_n\}$. Actually, the limit of the sequence of minimal immersions will be the complete proper minimal immersion we are looking for.

At this point, we are able to define the immersion that proves the theorem. This immersion will be the limit of the sequence $\{X_n \}_{n \in \n}$. But before of this, we need to construct the domain of definition of the limit immersion. 
Properties (I$_n$), $n \in \n$,  imply that the set:
$$\Omega=\bigcup_{n=1}^\infty \intc P_n^{\varepsilon_{n+1}}=\bigcup_{n=1}^\infty \intc P_n^{\xi_{n}}$$
is an expansive union of simply connected domains, resulting in $\Omega$ being simply connected. Moreover, as $\intc P_n^{\xi_{n}} \subset \d$, then $\Omega \subset \d$. Using Riemann's mapping Theorem, we deduce that $\Omega$ is conformally equivalent to the unit disk. 

Properties (II$_n$) say us that $\{X_n\}_{n \in \n}$ is a Cauchy sequence, uniformly on compact sets of $\Omega$, and so, by using Harnak's theorem, it converges. Let $\psi_{(\Gamma,\epsilon)}:\Omega\rightarrow \r^3$ be the limit of $\{X_n\}_{n \in \n}$. Then $\psi_{(\Gamma,\epsilon)}$ has the following properties:

\begin{claim} \label{cl:conformal}
$\psi_{(\Gamma,\epsilon)}$ is a conformal minimal immersion.
\end{claim}
\begin{proof} The proof of this claim is a direct consequence of Properties  (III$_n$), $n \in \n.$ \end{proof}

\begin{claim} \label{cl:propia}
 $\psi_{(\Gamma,\epsilon)}:\Omega \longrightarrow C$ is proper. 
\end{claim}
\begin{proof}
Consider a compact subset $K \subset C$. Let $n_0$ be a natural so that $$K \subset (E^{n-1})_{-2(b_{n-1}+b_n)-\sum_{k\geq n} \varepsilon_k}, \quad \forall n  \geq n_0.$$ From Properties (VI$_n$), we have $X_n(z)\in \r^3- (E^{n-1})_{-2(b_{n-1}+ b_n)}$, $\forall z \in \intc P_n- \intc P_{n-1}^{\varepsilon_n}$. Moreover, taking into account (II$_k$), for $k \geq n$, we obtain 
$$\psi_{(\Gamma,\epsilon)}(z)\in \r^3 - (E^{n-1})_{-2(b_{n-1}+ b_n)-\sum_{k\geq n} \varepsilon_k}.$$
Then, we have $\psi_{(\Gamma,\epsilon)}^{-1}(K)\cap(\intc P_n - \intc P_{n-1}^{\varepsilon_n})=\emptyset$ for $n\geq n_0$. This implies that $\psi_{(\Gamma,\epsilon)}^{-1}(K)\subset \intc P_{n_0-1}^{\varepsilon_{n_0}}$, and so it is compact in $\Omega$.
\end{proof}

\begin{claim} \label{cl:complete}
$\Omega$ is complete with the metric $\metri{ \psi_{(\Gamma,\epsilon)}}$. 
\end{claim}
\begin{proof}
This is a trivial consequence of Properties (II$_n$) and (IV$_n$), $n \in \n$. 
\end{proof}

\begin{claim} \label{cl:hausdorff} The limit set $\psi_{(\Gamma,\epsilon)}(\partial \Omega)$ satifies $\delta^H(\psi_{(\Gamma,\epsilon)}(\partial \Omega),\Gamma)<\epsilon$. 
\end{claim}
\begin{proof}
\noindent First, we are going to prove that $\dist_{\r^3}(p,\Gamma)< \tau$, for all $p$ in the limit set $\psi_{(\Gamma,\epsilon)}(\partial \Omega)$.

Given $p\in \psi_{(\Gamma,\epsilon)}(\partial \Omega)$, we know that there exists a sequence of points in $\Omega$, $\{z_n\}_{n \in \n}$, with 
$$\{ \psi_{(\Gamma,\epsilon)}(z_n)\}_{n \in \n} \to p.$$ 
Without loss of generality (up to re-index the sequence and take a subsequence) we can assume:
\begin{eqnarray}
z_n & \in & \intc P_n - \intc P_{n-1}^{\varepsilon_n} \label{eq:martima-1} \\
z_n & \in & \intc P_{n}^{\varepsilon_{n+1}}  \label{eq:martima-2}
\end{eqnarray}

From \eqref{eq:martima-1} and taking (VII$_n$) into account we obtain:
\begin{equation} \label{eq:esa-1}
\dist_{\r^3} \left(X_n(z_n),\Gamma \right) \leq \mu_n <\frac{\tau}2.
\end{equation}
On the other hand, using \eqref{eq:martima-2} and the fact that $z_n$ satisfies the conditions to apply Properties (II$_n$), for $k\geq n+1,$ then we have:
\begin{equation} \label{eq:esa-2}
\| \psi_{(\Gamma,\epsilon)}(z_n)-X_n(z_n)\| \leq \sum_{k=n+1}^\infty \|X_k(z_n)-X_{k-1}(z_n) \|< \frac{\tau}4.
\end{equation}
Combining \eqref{eq:esa-1} and \eqref{eq:esa-2} we trivially deduce that $\dist_{\r^3}(\psi_{(\Gamma,\epsilon)}(z_n), \Gamma) < 3\tau/4$. Finally, we take limit as $n \to \infty$ and obtain $\dist_{\r^3}(p, \Gamma)  \leq 3 \tau/4<\epsilon$ (recall that we have choosen $\tau < \epsilon.$)
\vskip .7cm

\noindent Now, we have to prove that $\dist_{\r^3}(x,\psi_{(\Gamma,\epsilon)}(\partial \Omega))<\tau,$ for all $x\in \Gamma$.

Pick a point $x \in \Gamma$. Property (VII$_n$) once again says to us:
$$\delta^H(X_n(\intc P_n -\intc P_{n-1}^{\varepsilon_n}), \Gamma) \leq \mu_n< \frac {\tau}2 .$$
Assume that $\b(x,\tau/2) \cap X_n(P_{n-1}^{\varepsilon_n}) $ is empty. It implies that $X_n(P_{n-1}^{\varepsilon_n})$ and $\Gamma$ can not be homologous in $T(\Gamma,\tau/2)$. But Properties (IX$_n$) and (VII$_n$) imply that $X_n(P_{n-1}^{\varepsilon_n})$ is homologous to $\Gamma$ in $T(\Gamma,\mu_n) \subset T(\Gamma,\tau/2)$, which is absurd.
 
 Hence, we deduce that there exists $z_n\in P_{n-1}^{\varepsilon_n}$ such that $\|X_n(z_n)-x\|< \tau/2.$ 
Moreover, reasoning as in \eqref{eq:esa-2}, we have $\| \psi_{(\Gamma,\epsilon)}(z_n)-X_n(z_n)\|<\tau/4$. Then, we conclude \begin{equation} \label{eq:ultima}
\|\psi_{(\Gamma,\epsilon)}(z_n)-x\|< 3 \tau/4, \quad \forall n \in \n.
\end{equation} 

Notice that the sequence $\{\psi_{(\Gamma,\epsilon)}(z_n)\}_{n \in \n}$ admits a convergent subsequence (recall that $\overline{C}$ is compact). Let $p$ denote the limit of such  subsequence. Inequality \eqref{eq:ultima} yields that $\|p-x\| \leq 3 \tau/4$. In particular, we have $\dist_{\r^3}(x, \psi_{(\Gamma,\epsilon)}(\partial \Omega))< \epsilon.$
\end{proof}

This completes the proof of Theorem \ref{th:mari}.
\end{proof}
\bigskip

As we mentioned in the introduction, the above theorem has the following direct application:

\begin{teorema} \label{th:regular}
Every bounded  domain with regular boundary admits a complete properly immersed minimal disk.
\end{teorema}
\begin{proof}
If  $R$ is a bounded regular domain, then there exists an open connected region $A \subset \partial R$ so that the mean and Gauss curvatures are positive on $A$. Consider $R'$ an strictly convex regular domain $R' \subset R$ and such that $\emptyset \neq \partial R' \cap \partial R\subset A$, like in Figure \ref{fig:regular}. Hence, in order to get the minimal  immersion we are looking for, it suffices to apply Theorem \ref{th:mari} to a Jordan curve $\Gamma$ in $\partial R' \cap \partial R$ and an $\varepsilon >0$ small enough.
\begin{figure}[hb]
	\begin{center}
		\includegraphics[width=0.46\textwidth]{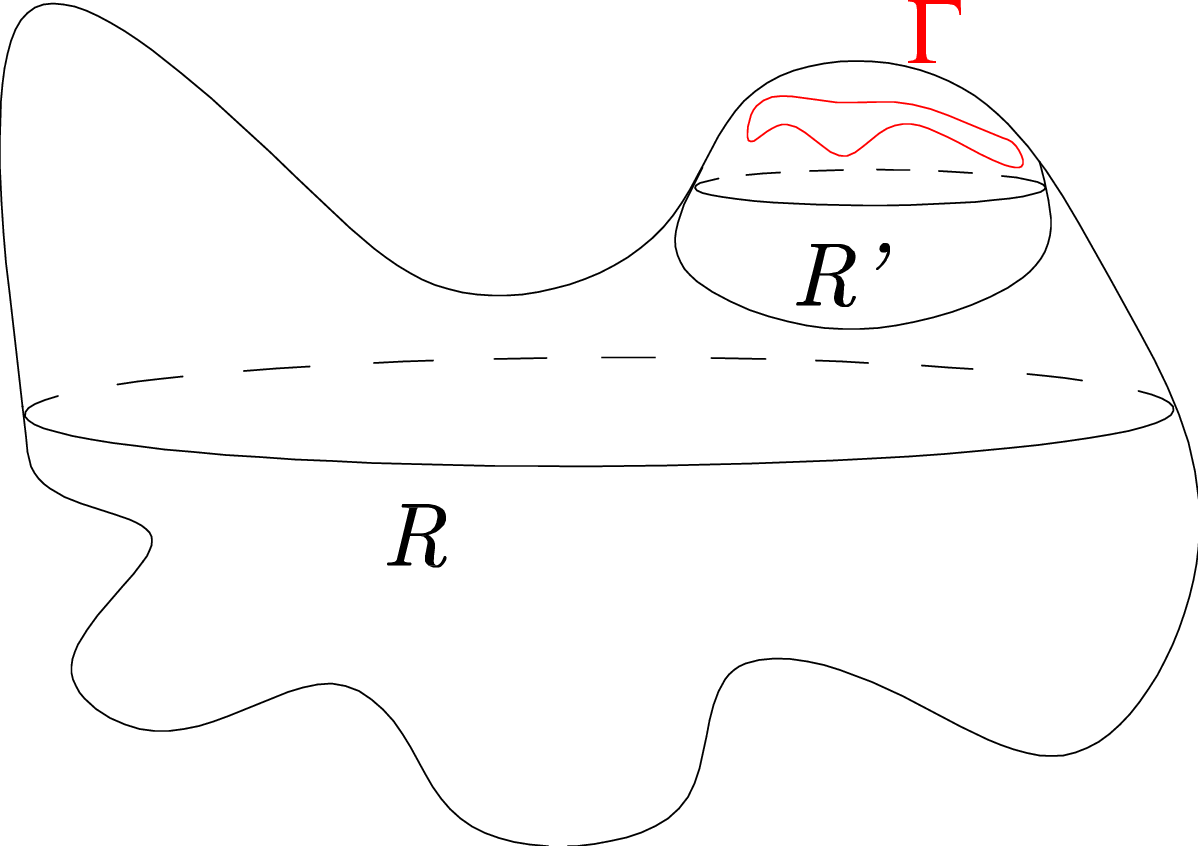}
	\end{center}
	\caption{The domains $R$ and $R'$ and the curve $\Gamma.$}
	\label{fig:regular}
\end{figure}

\end{proof}
Theorem \ref{th:regular} implies that a bounded regular domain is not {\em universal for minimal surfaces.} A connected region of space which is open or the closure of an open set is {\it universal for minimal surfaces}, if every complete properly immersed minimal surface in the region is recurrent for Brownian motions. In particular, a bounded domain is universal if and only if it contains no complete properly immersed minimal surfaces.

The other interesting consequence follows from a well known result in convex geometry.
 Next theorem essentially asserts that the Jordan curve $\Gamma$ can be substituted by an arbitrary compact set in the statement of Theorem \ref{th:mari}.

\begin{teorema} \label{th:compacto}
Let $C$ be an  regular, strictly convex  bounded domain, and consider a {\em connected compact set} $K \subset \partial C$. Then, for each $\varepsilon>0$, there exists a complete proper minimal immersion $\varphi_{(K,\varepsilon)}:\d \rightarrow C$ satisfying that the Hausdorff distance $\delta^H(\psi_{(K,\varepsilon)}(\partial \d), K) < \varepsilon.$
\end{teorema}
Theorem \ref{th:compacto} is a trivial consequence of Theorem \ref{th:mari} and the lemma  below.
\begin{lema}
Let $K$ be a connected compact set in $\partial C$, then for every $\nu>0$ there exists an smooth Jordan curve $\Gamma$ so that $\delta^H(K,\Gamma)< \nu.$
\end{lema}
\begin{proof}
Classical results about Hausdorff metric says to us that there is a finite set $F \subset K$ such that $\delta^H(K,F) < \nu/3.$ Label $K'=\left(K+\b(0,\nu/3)\right) \cap \partial C$, which is a connected open set in $\partial C.$ Then, it is clear that we can find a compact simple piecewise smooth curve $\gamma'$ in $K'$ passing through all the points in $F.$ If $\gamma'$ is not closed, we can consider a curve parallel to $\gamma'$, that we call $\gamma''$,  and sufficiently close to $\gamma'$ so that the Jordan curve $\gamma$ obtained joining the extremes of $\gamma'$ and $\gamma''$ satisfies $\delta^H(\gamma, K)<\nu/3.$ To finish, we only need to approximate $\gamma$ by a smooth Jordan curve satisfying $\delta^H(\gamma,\Gamma)<\nu/3.$
\end{proof}

Finally, we would like to show that our results are sharp in the following sense. If we remove from Theorem \ref{th:mari} the hypothesis of $C$ being strictly convex, then the result fails. The most simple counterexample is open halfspace. If we consider a complete  minimal disk properly immersed in an open halfspace then the limit set cannot be bounded, if not we will arrive to a contradiction by the maximum principle. Actually we can prove that:
\begin{proposicion} \label{pro:fails}
Let $D$ a domain of $\rth$ satisfying that there exists an open region $U \subset \partial D$ where the mean curvature associated to the inward pointing normal is non positive. If $U$ is a graph over a bounded convex domain of a plane, then there are no complete proper minimal disks in $D$ whose limit set lies on $U$. In particular, Theorem \ref{th:mari} is not true for domains of this kind.
\end{proposicion} 
\begin{proof}
The proof of this proposition is an easy application of the maximum principle. We proceed by contradiction. Assume there exists a complete proper minimal immersion $\psi:\d \rightarrow D$ so that the limit set $\psi(\partial \d)$ is contained in $U.$ From the hypotheses, we know that $U$ is a graph over a plane $\Pi.$ First, we translate $U$, orthogonally to $\Pi$, and toward the interior of $D$ until it does not touch $\psi(\d)$ (see Figure \ref{fig:maximo}.)
\begin{figure}[htbp]
	\begin{center}
		\includegraphics[width=0.46\textwidth]{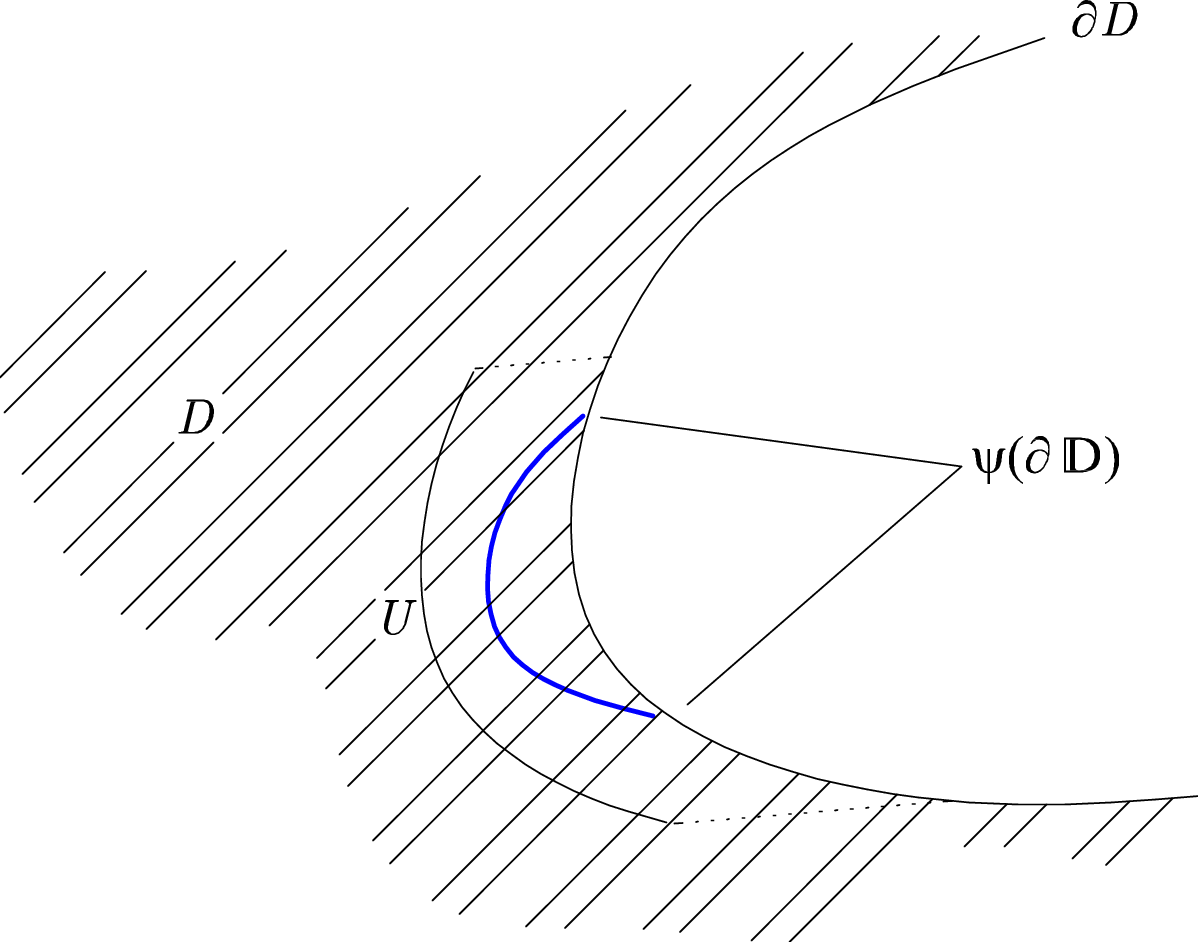}
	\end{center}
	\caption{We translate $U$ in the direction of the interior of $D$  until it does not touch $\psi(\d)$.}
	\label{fig:maximo}
\end{figure}
Now continuously translate $U$ toward $\psi(\d)$ until the translated graph intersects $\psi(\d)$ for the first time. This must occur because the immersion is proper. As the limit set is contained in $U$, then $\psi(\d)$ is contained in the interior of the convex cylinder over the projection of $U$ in $\Pi$. In particular the boundary of the graph never touches $\psi(\d)$, then the translated graph $U_0$ and  $\psi(\d)$ have an interior contact point, and so the maximum principle leads to $U_0=\psi(\d)$ which is absurd. This contradiction  proves the proposition.
\end{proof}

On the other hand, regularity is also necessary. Indeed, Martín, Meeks and Nadirashvili \cite{mmn} have recently proved that:
\begin{teorema}[\cite{mmn}] \label{th:first}
Let ${\cal D}$ be any bounded open domain in $\r^3$. Then there exists a
proper countable collection $\cal F$ of pairwise disjoint horizontal
simple closed curves in $\cal D$ such that the complementary domain
 $\widetilde{\cal D}={\cal D} - {\cal F}$ is universal for minimal surfaces
  with at least one annular end. In particular, any complete immersed minimal
   surface of finite genus in $\widetilde{\cal D}$ must have an uncountable number of ends.
\end{teorema}
Among other things, this means that Theorem \ref{th:regular} is sharp, too.

 \vspace{10mm}

\noindent Francisco Martín, Santiago Morales \\ Departamento de Geometría y Topología \\  Universidad de Granada \\ 18071 Granada, Spain  \\
\texttt{fmartin@ugr.es, santimo@ugr.es}
\end{document}